\magnification=\magstep1

\catcode`@=11
\def\displaylines#1{\displ@y
  \halign{\hbox to\displaywidth{$\@lign\hfil\displaystyle{}##\hfil$}\crcr
    #1\crcr}}
\catcode`@=12

\newcount\contaparag
\newcount\contaeq

\let\TeXeqno\eqno
\let\TeXeqalignno\eqalignno

\def\eqno#1{\global\advance\contaeq 1
  \ifx#1\undefined
  \xdef#1{\bgroup\noexpand\rm
    (\the\contaparag.\the\contaeq)\egroup}%
  \wlog{\def\string#1{(\the\contaparag.\the\contaeq)}}%
  \else
    \def\gianota{#1}%
    \ifx\gianota\space
    \else
      \messaggiogianota
    \fi
  \fi
  \TeXeqno#1}

\def\messaggiogianota{\wvd{}\wvd{ONE HAS ALREADY GIVEN A MEANING TO
  \expandafter\string\gianota\space }\wvd{}}

\def\wvd#1{\immediate\write 16 {#1}}

\let\TeXeqalignno\eqalignno

\def\eqalignno#1{\defs#1\finedefs{#1}}

\def\defs#1&#2&#3\cr#4{%
  \def\arg{#3}%
  \ifx#3\undefined
    \global\advance\contaeq 1
    \xdef#3{\bgroup\noexpand\rm
      (\the\contaparag.\the\contaeq)\egroup}%
    \wlog{\def\string#3{(\the\contaparag.\the\contaeq)}}%
  \else
    \def\gianota{#3}%
    \ifx\gianota\space
    \else
      \messaggiogianota
    \fi
  \fi
  \ifx#4\finedefs
    \let\next\finedefs
  \else
    \def\next{\defs#4}%
  \fi
  \next}

\let\finedefs\TeXeqalignno

\def\accorpa #1#2{%
  {%
   \rm\expandafter\expandafter\expandafter\ccrp\expandafter #1#2}}

\def\ccrp\bgroup\rm(#1.#2)\egroup
  \bgroup\rm(#3.#4)\egroup{(#1.#2-#4)}

% A ogni paragrafo nuovo usare \nuovoparagrafo seguito dal numero del
% paragrafo prima di cio' che c'e' gia' scritto (che andra' a capo
% senza righe vuote intermedie)

\def\nuovoparagrafo #1 {\contaparag #1 \contaeq 0 }

% numerazione dei riferimenti bibliografici

\newcount\contaref

\def\numeraref#1{%
  \def\argref{#1}%
  \ifx\argref\fineelenco
    \let\next\relax
  \else
    \advance\contaref 1
    \edef#1{\the\contaref}%
    \wlog{\def\string#1{#1}}%
    \let\next\numeraref
  \fi
  \next}

\def\fineelenco{\fineelenco}

%%%%%%%%%%%%%%%%%%%%%%%%%%%%%%%%%%%%%%%%%%%%%%%%%%%%%%%%%%%%%

\mathsurround=2pt
\def\wbox#1/{\vbox{\hrule\hbox{\vrule height#1\kern#1\vrule height#1}\hrule}}
\newbox\boxfine \setbox\boxfine=\hbox{\wbox6pt/}
\def\fine{\ \kern6pt\copy\boxfine}
\mathchardef\TeXchi="011F
\font\maiuscolo=cmcsc10
\font\piccolo=cmr10 at 8 pt
\font\piccolorm=cmr10 at 7pt
\font\piccoloit=cmti10 at 7 pt
\font\tito=cmbx10 scaled \magstep2
\font\titosl=cmbxsl10 scaled \magstep2
\font\tenmsb=msbm10  
\font\macc=cmtt10
\font\sevenmsb=msbm7
\font\fivemsb=msbm5
\newfam\msbfam
\textfont\msbfam=\tenmsb
\scriptfont\msbfam=\sevenmsb
\scriptscriptfont\msbfam=\fivemsb
\def\Bbb{\fam\msbfam }
 \def\enne{{\Bbb N}}  
\def\erre{{\Bbb R}}
\def\div{{\rm div}}
\def\chi{{\setbox0 =\hbox{$\mathsurround=0pt\TeXchi$}\hbox
  {\raise\dp0 \copy0 }}}

\def\dual#1{\langle#1\rangle}

\def\teta{\vartheta}

\def\eps{\varepsilon}

\def\mez{{1\over 2}}
\def\quarto{{1\over 4}}

\def\pan{{\cal P}^n_{g_0}}
\def\diesis{2^{\#}}
\def\fix{{\rm Fix\,}}
\def\laplconf{L_{g_0}}
\def\pii{p}

\def\ui{u_i}
\def\yi{y_i}
\numeraref
\pertuno
\pertdue
\spagnoli
\alm
\simmetria
\aubin
\bgm
\branson
\bcy
\chang
\cgydue
\cgy
\cqy
\changyangdue
\changyang
\dhl
\dma
\dmadue
\gursky
\hebey
\liuno
\lidue
\lincs
\paneitz
\van
\weixu
\fineelenco

\nopagenumbers
\quad
\vskip2truecm\noindent
\nuovoparagrafo 0
\centerline{\tito Existence of conformal metrics on {\titosl S}$^{{\displaystyle^n}}$}
\medskip\noindent
\centerline{\tito with prescribed fourth order invariant}
\vskip1truecm
\centerline{\bf Veronica Felli\footnote{$^{(1)}$}{\piccolo Supported by M.U.R.S.T. under the national project ``Variational Methods and Nonlinear Differential Equations.''.}}
\vskip1truecm
\centerline{S.I.S.S.A.}
\centerline{via Beirut, 2-4}
\centerline{34014 Trieste, Italy}
\centerline{e-mail: {\macc felli@sissa.it}} 
\vskip2truecm\noindent
\centerline{\vbox{\hsize=14truecm{\bf \noindent Abstract.}\quad{\piccolorm In this paper we prescribe a fourth order conformal invariant on the standard {\piccoloit n}-sphere, with ${\scriptstyle n\geq5}$, and study the related fourth order elliptic equation. We first find some existence results in the perturbative case. After some blow up analysis we build a homotopy to pass from the perturbative case to the non-perturbative one under some flatness condition. Finally we state some existence results under the assumption of symmetry.}}}
\footline{\hfil\folio\hfil}
\vskip2truecm\noindent
\centerline{\bf 0.\quad Introduction.}
\bigskip\noindent
Let $(M^4,g)$ be a smooth 4-dimensional manifold, $S_g$ the scalar curvature of $g$, ${\rm Ric}_g$ the Ricci curvature of $g$ and $\Delta_g=\div_g\,\nabla_g$ the Laplace-Beltrami operator on $M^4$. Let us consider the fourth order operator discovered by Paneitz [\paneitz] in 1983:
$${\cal P}^4_g\psi=\Delta^2_g\psi-\div_g\,\left({2\over 3}S_g-2\,{\rm Ric}_g\right)d\psi.$$
The Paneitz operator is conformally invariant on 4-manifolds, that is if $g_\omega=e^{2\omega}g$ then
$${\cal P}_{g_\omega}(\psi)=e^{-4\omega}{\cal P}_g(\psi)\quad\hbox{for all}\quad \psi\in C^\infty(M^4)$$
and it can be considered as a natural extension of the conformal laplacian on 2-manifolds. It is also known that 
$${\cal P}_g^4\omega+2Q_{g_\omega}e^{4\omega}=2Q_g$$
where
$$Q_g={1\over 12}\left(-\Delta_gS_g+S_g^2-3|{\rm Ric}_g|^2\right).$$
Such a $Q_g$ is a fourth order invariant called {\sl $Q-$curvature}, because is the analogous for the Paneitz operator of the scalar curvature. It is natural to pose the problem of prescribing  the $Q-$curvature on $S^4$: given a smooth function $\widehat Q$ on $S^4$, to find a conformal metric $g_\omega$ such that $Q_{g_\omega}$ is $\widehat Q$. This problem leads to the following equation
$$ {\cal P}_{g_0}^4\omega+2\widehat Qe^{4\omega}=2Q_{g_0}\quad\hbox{on}\ (S^4,g_0)$$
which has been studied in [\weixu]. For more details about the basic properties of the Paneitz operator one can see [\chang] and [\changyang]. We also refer to [\bcy, \cgy, \cqy, \gursky] for results about ${\cal P}^4_{g_0}$.\par
The Paneitz operator was generalized to higher dimension by Branson [\branson]. Given a smooth compact Riemannian $n-$manifold $(M^n,g)$, $n\geq 5$, the Branson-Paneitz operator is defined as
$${\cal P}^n_g\psi=\Delta^2_g \psi-\div_g\,(a_n S_g g+b_n {\rm Ric}_g)d\psi+{n-4\over 2}Q_g^n\psi,$$
where
$$\displaylines{a_n={(n-2)^2+4\over 2(n-1)(n-2)},\qquad b_n=-{4\over n-2},\cr
Q_g^n=-{1\over 2(n-1)}\Delta_gS_g+{n^3-4n^2+16n-16\over 8(n-1)^2(n-2)^2}S_g^2-{2\over (n-2)^2}|{\rm Ric}_g|^2.}$$
See [\dhl] for details about the properties of ${\cal P}_g^n$. If $g_v=v^{4\over n-4}g$ is a conformal metric to $g$, one has that for all $\psi\in C^\infty(M)$,
$${\cal P}^n_g(\psi v)=v^{n+4\over n-4}{\cal P}^n_{g_v}\psi$$
and
$${\cal P}^n_gv={n-4\over 2}Q^n_{g_v}v^{n+4\over n-4}. \eqno\zerouno$$
\footline{\hfill\folio\hfill}
It is natural to study the problem of prescribing $Q$ on the standard sphere $(S^n,g_0)$, that is of finding solutions to equation \zerouno\ with $Q^n_{g_v}$ equal to some prescribed function $\widehat Q$, i.e.
 $${\cal P}^n_{g_0}v={n-4\over 2}\widehat Q v^{n+4\over n-4},\quad v>0.\eqno(P) $$
Problem $(P)$ can be viewed as the analogue of the classic scalar curvature problem on $(S^n, g_0)$. While the scalar curvature problem has been widely studied, (see for example the monographes [\aubin, \hebey], see also [\liuno, \lidue]), problem $(P)$ concerning the $Q-$curvature has been faced in [\dhl, \dma, \dmadue] for dimension $n=5,6$ only. The purpose of the present paper is to study $(P)$ and prove results that are the counterpart of those known for the scalar curvature problem.\par
In Section 1 we start proving an existence result for equation $(P)$ on $(S^n,g_0)$ in the case $\widehat Q=1+\eps K$ (see Section 1, Theorem 1.5) using a perturbative method introduced in [\pertuno] and [\pertdue] and generalizing to the Branson-Paneitz operator the results of [\spagnoli]. This result is a slight improvement of [\dma, Theorem 1.3], where problem $(P)$ is solved in the perturbative case with some more assumptions on $\widehat Q$.
\par
In Section 2 we extend to $(P)$ the results of [\liuno]. We first adapt the arguments of [\liuno] to carry out some blow up analysis for a family of equations approximating $(P)$. In particular, in Theorem 2.12 we prove that, making some flatness assumption on $\widehat Q$, no blow up occurs and hence any positive solution to problem $(P)$ with $\widehat Q$ replaced by $\widehat Q_t=t\widehat Q+(1-t){(n^2-4)n/ 8}$ stays in a compact set of $C^{4,\alpha}(S^n)$, $0<\alpha<1$. This compactness result allows us to pass from the perturbative case treated in Section 1 to the non-perturbative one by means of a homotopy argument like in [\liuno]. In this way we obtain an existence result for $(P)$ (Theorem 2.13). Let us point out that in [\dmadue] no flatness condition is required and blow up can occur when $n=6$. A more precise comparison with the results of [\dmadue] is made in Remark 2.14.  
\par
Finally, in Section 3 we consider the symmetric case and generalize to the Branson-Paneitz operator the results of [\simmetria] in the perturbative case and of [\alm] in the non-perturbative one.
\par Our techniques are strongly based on those of [\dmadue] and [\liuno]. As far as [\liuno] is concerned, we assume that the reader is familiar with the results and the techniques of this work, while we recall the formulas we need from [\dmadue].
\bigskip\noindent
\centerline{\bf Acknowledgements}
\medskip
The author wishes to thank Prof. A. Ambrosetti for having proposed her the study of this problem and for his useful advices, and Dr. M. Ould Ahmedou and Dr. A. Malchiodi for helpful suggestions.

\vskip2truecm
\nuovoparagrafo 1

\centerline{\bf 1.\quad Perturbation existence results}
\bigskip\noindent
We consider the Paneitz operator on the standard unit sphere $(S^n,g_0)$
$$ {\cal P}^n_{g_0}v=\Delta_{g_0}^2v-c_n\Delta_{g_0}v+d_nv \eqno \uno$$ where
$n\geq 5$, $\Delta_{g_0}$ stands for the Laplace-Beltrami operator on $S^n$, and$$c_n=\mez(n^2-2n-4)\qquad d_n={n-4\over 16}n(n^2-4).$$
We deal with problem $(P)$ on $S^n$. Denoting the euclidean metric
on $\erre^n$ by $\delta$ and the stereographic projection on $\erre^n$ by $\pi$, one has 
$$(\pi^{-1})^{\star} g_0=\varphi^{4\over n-4}\delta$$ where 
$$\varphi(y)=\left({2\over 1+|y|^2}\right)^{n-4\over 2}.\eqno \dueprimo$$ Using this conformal change of metric and the conformal covariant properties of ${\cal P}^n$ one has
$$\pan (\Phi(u))=\varphi^{-{n+4\over n-4}}{\cal
P}^n_\delta(u)=\varphi^{-{n+4\over n-4}}\Delta^2 u\qquad\forall\, u\in
C^\infty(\erre^n), \eqno\conf $$
where
$$\Phi\ :\quad{\cal D}^{2,2}(\erre^n)\ \longrightarrow\ H^2(S^n),\qquad
                     \Phi(u)(x)={u(\pi(x))\over \varphi(\pi(x))}$$
is an isomorphism between $H^2(S^n)$ and $E:={\cal D}^{2,2}(\erre^n)$, the completion of
$C^\infty_c(\erre^n)$ w.r.t. the norm 
$\|u\|^2=\int_{\erre^n}|\Delta u|^2$.
\par\noindent
After the transformation $v=\Phi(u)$ and letting $Q=\widehat Q\circ \pi^{-1}$, problem $(P)$ becomes, up to an
uninfluent constant,
 $$\left\{\eqalign{&\Delta^2u= Q u^{n+4\over n-4}\cr
                   &u>0,\qquad\hbox{in}\quad \erre^n.\cr}\right.$$
We take $Q$ of the form $1+\eps K$ for some bounded function $K$, such that problem $(P)$ in the perturbative case after
stereographic projection is the following
$$\Delta^2 u=(1+\eps K)u^p,\qquad u>0\quad\hbox{in}\ \erre^n\eqno\tre$$
where $$p=2^{\#}-1={n+4\over n-4}\quad\hbox{and}\quad 2^{\#}={2n\over n-4}.$$
We recall that by the Sobolev embedding theorem, we have the embedding
$H^2(S^n)\hookrightarrow L^{\diesis}(S^n)$ which is not compact. Consider the
functional $f_\eps:E\to\erre$ given by
$$ f_\eps(u)=\mez \int_{\erre^n}(\Delta u)^2 -{1\over
p+1}\int_{\erre^n}|u|^{p+1}-{\eps\over p+1}\int_{\erre^n}K(x)|u|^{p+1}.$$
Plainly, $f_\eps\in C^2(E,\erre)$ and any critical point $u\in E$ is a solution of
$$\Delta^2u=(1+\eps K)|u|^{p-1}u\quad\hbox{in}\ \erre^n.\eqno \quattro$$
\medskip\noindent
{\bf 1.1. Critical points of $f_\eps$.}\quad It is known (see [\lincs ]) that the positive solutions of the unperturbed problem 
$$\Delta^2u=|u|^{p-1}u$$
 are given by the functions
$$z_{\mu,\xi}(y)=\mu^{-{n-4\over2}}z_0\left({y-\xi\over\mu}\right),\qquad \mu\in\erre^+,\ \xi\in\erre^n$$
where
$$z_0(y)=C_n(1+|y|^2)^{-{n-4\over2}},\qquad C_n=[n(n^2-4)(n-4)]^{n-4\over8}.$$
We set 
$$Z=\{z_{\mu,\xi}:\ \mu>0,\ \xi\in \erre^n\}.$$
We can face problem \tre\ with the same approach used in [\spagnoli] for the scalar curvature problem, i.e. the abstract perturbation method discussed in [\pertuno] and [\pertdue].\par
In order to apply the abstract setting of [\pertuno] and [\pertdue], we need the following result.
\medskip\noindent
{\bf Lemma 1.1.} \quad \sl $f_0$ satisfies the following properties
\item{\rm(i)}\quad $\dim Z=n+1$;
\item{\rm(ii)}\quad $D^2f_0(z)=I-C\quad\forall\, z\in Z$, where $I$ is the identity and $C$ is compact;
\item{\rm(iii)}\quad $T_{z_{\mu,\xi}}Z=\ker\{D^2f_0(z_{\mu,\xi})\}\quad\forall \,\mu>0,\ \xi\in\erre^n$.
\medskip\noindent\rm
{\bf Proof.}\quad It is sufficient to give only the proof of (iii), which is a {\sl nondegeneracy condition}. Since the inclusion $T_{z_{\mu,\xi}}Z\subseteq\ker\{D^2f_0(z_{\mu,\xi})\}$ is obvious, it is enough to prove that
$$\dim T_{z_{\mu,\xi}}Z=\dim \ker\{D^2f_0(z_{\mu,\xi})\}.$$
Up to a translation, we can assume that $\xi=0$ and for simplicity, we consider $\mu=1$. We have
$$u\in\ker \{D^2f_0(z_0)\}\Longleftrightarrow {\cal P}^n_{\delta}u=\Delta^2u=pz_0^{p-1}u\quad\hbox{in}\ \erre^n$$
that is, using \conf, if and only if
$$\varphi^{n+4\over n-4}\pan(\Phi(u))={n+4\over n-4}z_0^{p-1}u.$$
Hence $u\in\ker\{D^2f_0(z_0)\}$ if and only if $v=\Phi(u)$ satisfies
$$\pan v={n+4\over n-4}z_0^{p-1}\varphi^{-{8\over n-4}}v={1\over 16}n(n^2-4)(n+4)v\quad\hbox{in}\ S^n.\eqno\nondeg$$
Recalling \uno, equation \nondeg\ can be written as
$$\Delta^2_{g_0}v-c_n\Delta_{g_0}v=\left[{1\over 16}n(n^2-4)(n+4)-{1\over 16}n(n^2-4)(n-4)\right]v=\mez n(n^2-4)v.$$
Therefore
$\ker\{D^2f_0(z_0)\}=\Phi^{-1}(V)$
where $V$ is the eigenspace of $\Delta^2_{g_0}-c_n\Delta_{g_0}$ corresponding to the eigenvalue $\mez n(n^2-4)$. It is easy to check \footnote{$^{(2)}$}{\piccolo An explicit proof of this is given in the preprint of [\dhl]; in the final version it has been eliminated.} that there is a correspondence between the spectrum of the Paneitz operator and the spectrum of the Laplace-Beltrami operator: for any $\alpha>0$ one has
$${\rm Sp}\,(\Delta^2_{{g_0}}-\alpha\Delta_{{g_0}})=\{\lambda^2+\alpha\lambda:\ \lambda\in {\rm Sp}\,(-\Delta_{g_0})\}.$$
Moreover $v$ is an eigenfunction of $\Delta^2_{{g_0}}-\alpha\Delta_{{g_0}}$ associated to $\lambda^2+\alpha\lambda$ if and only if $v$ is an eigenfunction of $-\Delta_{g_0}$ associated to $\lambda$.\par\noindent
In our case, we take $\alpha=c_n$ and hence from $\lambda^2+\alpha\lambda=\mez n(n^2-4)$ we infer $\lambda=n$. Then $V$ concides with the eigenspace of $-\Delta_{g_0}$ associated to $n$. It is well-known [\bgm] that the spectrum of the Laplace-Beltrami operator on the $n-$standard sphere $S^n$ is given by
$$\lambda_k=k(n+k-1),\qquad k\geq0$$
and the dimension of the eigenspace associated to $\lambda_k$ is
$${(n+k-2)!(n+2k-1)\over k!(n-1)!}\qquad\forall\, k\geq 0.$$
\rm Since $\lambda_1=n$, one has that $\dim V=n+1$
hence $\dim \ker \{D^2f_0(z_0)\}=n+1=\dim T_{z_0}Z$ and this gives (iii).\fine\medskip\noindent 
Now we follow closely the procedure of [\spagnoli] and assume that the reader is quite familiar with it. The previous lemma allows us to apply the abstract method and to define, for $\eps$ small, a smooth function $w_\eps(z):\ Z\to(T_zT)^{\bot}$ such that the perturbed manifold $Z_\eps=\{z_{\mu,\xi}+w_\eps(z_{\mu,\xi})\}$ is a natural constraint for $f_\eps$. There results
$$\varphi_\eps(z_{\mu,\xi}):=f_\eps(z_{\mu,\xi}+w_\eps(z_{\mu,\xi}))=b-\eps\Gamma(\mu,\xi)+o(\eps)\qquad\eps\to 0$$
where 
$$ b:={2\over n}\int_{\erre^n}z_0^{\diesis}=f_0(z)\qquad\forall\, z\in Z$$
and
$$\Gamma(\mu,\xi)={1\over\diesis}\int_{\erre^n}K(y)z_{\mu,\xi}^{\diesis}(y)\,dy={1\over p+1}\int_{\erre^n}K(\mu\zeta+\xi)z_0^{p+1}(\zeta)\, d\zeta.$$
Thanks to [\spagnoli, Theorem 2.1 and Remark 2.2], if $\Gamma$ satisfies
$$\leqalignno{&\Gamma \ \hbox{has an isolated set of critical points}\ {\cal C}\ \hbox{such that}\ \deg(\Gamma',\Omega,0)\not=0, &(A)\cr
& \hbox{where}\ \Omega\ \hbox {is an open bounded neighbourhood of}\  {\cal C},&\cr}$$
then for $\eps$ small enough, the functional $f_\eps$ has a critical point $u_\eps$ such that $u_\eps$ approaches the critical manifold $Z$ as $\eps \to 0$. Since $\Gamma$ has the same form as in [\spagnoli], we can repeat one of the arguments of Section 3 of [\spagnoli] to handle with one of the possible cases which have been treated there. Here we state just one result for the sake of brevity and remark that also the other results of [\spagnoli] can be generalized to our problem.   
\medskip\noindent
{\bf Lemma 1.2.}\quad \sl Assume that $K$ satisfies
$$\left\{\eqalign{&\hbox{(K1.a)}\quad\exists\, \rho>0:\ \dual{K'(y),y}<0\quad\forall\,|y|\geq\rho,\cr
 &\hbox{(K1.b)}\quad\dual{K'(\cdot),\cdot}\in L^1(\erre^n)\quad\hbox{and}\quad \int_{\erre^n}\dual{K'(y),y}\, dy<0,\cr
&\hbox{(K1.c)}\quad K\in L^{\infty}(\erre^n)\cap C^1(\erre^n)\cr
           &\hbox{(K1.d)}\quad K\ \hbox{has finitely many critical points}\cr}\right.\leqno {(K1)}$$
and that
$$\leqalignno{&\forall\,\xi\in {\rm Crit}\, (K)\ \exists\,\beta\in]1,n[ \ \hbox{and}\  a_j\in C(\erre^n) \ \hbox{with}\  \tilde A_\xi:=\sum a_j(\xi)\not=0 \ \hbox{and}&(K2)\cr
&\hbox{such that}\ 
K(y)=K(\eta)+\sum a_j|y-\eta|^\beta+o(|y-\eta|^\beta)\ \hbox{as}\  y\to\eta\ \ \ \ \forall\,\eta&\cr
& \hbox{locally near}\  \xi. \ \hbox{Moreover there results}&\cr}$$ 
$$\sum_{\tilde A_\xi<0}\deg_{{\rm loc}}(K',\xi)\not=(-1)^n.$$
Then for any $\eps>0$ small, $f_\eps$ has a critical point $u_\eps\in{\cal D}^{2,2}(\erre^n)$ near $Z$; in particular $u_\eps$ is a solution of \quattro.
\medskip\noindent\rm
{\bf 1.2. Positivity.}\quad Now, following the techniques by Van Der Vorst [\van] as in [\dma], we prove that $u_\eps$ is strictly positive for $\eps$ sufficiently small, so that it is solution of problem \tre.
\par In order to prove $u_\eps>0$ we come back to $S^n$ and write problem \quattro\ in the form
$$\pan v_\eps=(1+\eps \widehat K)|v_\eps|^{\diesis-2}v_\eps,\eqno\posi$$
where $v_\eps=\Phi(u_\eps)$ and $\widehat K=K\circ \pi$. We recall that the volume element with respect to $g_0$ is
$$dV_{g_0}=|g_0|^{\mez}\,dy=\left[{2^{n-4\over 2}\over C_n}\right]^{2n\over n-4} z_0^{\diesis} \,dy$$
and observe that the solutions $u_\eps$ of \quattro\ given by the previous lemma approach some $\bar z\in {\rm Crit}\,(\Gamma)$ or some isolated set of critical points of $\Gamma$ included in $Z$. For simplicity, let us assume that $u_\eps$ converges to $z_0$ in $E$, hence $v_\eps=\Phi(u_\eps)$ converges in $H^2(S^n)$ to some positive constant $\bar c_n$ satisfying $\bar c_n^{\diesis-2}=d_n$. In particular, thanks to the Sobolev inequality, $v_\eps$ converges to $\bar c_n$ in $L^{\diesis}(S^n)$ and in measure. 
\bigskip\noindent{\bf Lemma 1.3.}\quad\sl $v_\eps$ is bounded in $L^\infty(S^n)$.\rm
\bigskip\noindent{\bf Proof.}\quad
Let us set
$$\Omega_\eps=\{x\in S^n:\ |v_\eps(x)|<2\bar c_n\}.$$
Since  $v_\eps$ converges to $\bar c_n$ in measure, the limit of $|S^n\,\backslash\,\Omega_\eps|$ as $\eps$ goes to $0$ is $0$. Now, after noting that
$$\pan v=L_{g_0}^2v-v$$
where $\laplconf v=-\Delta_{g_0}v+{c_n\over 2}v$, we can write \posi\ as
$$\laplconf^2 v_\eps=q_\eps|v_\eps|^{\diesis -2}v_\eps+g_\eps$$
where
$$\eqalign{&q_\eps(x)=\left(1+\eps \widehat K(x)+{1\over |v_\eps(x)|^{\diesis-2}}\right)\chi_{S^n\,\backslash\,\Omega_{\eps}}(x)\cr
           &g_\eps(x)=[(1+\eps \widehat K(x))|v_\eps(x)|^{\diesis -2}+1]v_\eps(x)\chi_{\Omega_\eps}(x).\cr}$$
Remark that, since on $S^n\,\backslash\,\Omega_{\eps}$ $v_\eps$ is far away from 0,  $q_\eps$ is bounded in $L^\infty$. Moreover one has
$$\|q_\eps|v_\eps|^{\diesis -2}\|^{n\over 4}_{L^{n\over 4}}=\int_{S^n\,\backslash\,\Omega_{\eps}}\left(1+\eps \widehat K+{1\over |v_\eps|^{\diesis-2}}\right)^{n\over 4}|v_\eps|^{\diesis}\, dV_{g_0}\ \mathop{\longrightarrow}\limits_{\eps\to 0} \ 0\eqno\posuno$$
and
$$\|g_\eps\|_{L^{\infty}}<\hbox{const}<\infty\quad\forall\,\eps>0.\eqno\posdue$$
We now use the following proposition which is proved in [\dma] using the techniques of [\van].
\medskip\noindent
{\bf Proposition 1.4.} \quad \sl Let $p\in [1, \diesis-1]$, $q\in L^\infty(S^n)$, and $g\in L^\infty(S^n)$. Suppose that $v\in H^2(S^n)$ is a weak solution of the equation
$$\laplconf^2 v=q|v|^{p-1}v+g,\qquad \hbox{in}\ S^n.$$
Then, for any $s>pn/4$, there exists a positive constant $\delta_n$ depending only on $n$, and a positive constant $C$ depending on $n$, $s$, and $\|q\|_{L^{\infty}}$, with the following properties. If $q|v|^{p-1}\in L^{n\over 4}(S^n)$ and $ \|q|v|^{p-1}\|_{L^{n\over 4}}\leq \delta _n$, then $v\in L^\infty(S^n)$ and its $L^\infty$ norm can be estimated as 
$$ \|v\|_{L^\infty}\leq C(\|g\|_{L^s}+\|g\|^p_{L^s}).\fine$$
{\bf Proof of Lemma 1.3 completed.}\quad \rm Taking $p=\diesis-1$, in view of \posuno\ and \posdue\ we can apply the previous proposition and deduce that $\{v_\eps\}_{\eps>0}$ is bounded in $L^{\infty}(S^n)$ with a positive constant not depending on $\eps$. \fine\medskip\noindent
Let us write \posi\ in the form
$$\laplconf^2v_{\eps}=b_\eps v_\eps,\eqno\ellegizero$$
where
$$b_\eps=(1+\eps \widehat K)|v_\eps|^{\diesis-2}+1,$$
or in the equivalent form
$$\laplconf^2(v_\eps-\bar c_n)=b_\eps(v_\eps-\bar c_n)+\bar c_n\left(b_\eps-{c_n^2\over 4}\right).\eqno \ellegizeroprimo$$
Owing to the convergence of $v_\eps$ to $\bar c_n$ and the boundness of $\{v_\eps\}$ in $L^{\infty}(S^n)$, by  dominated convergence there results that $b_\eps\to {\bar c_n}^{\diesis-2}+1={c_n^2\over 4}$ in $L^s$ for any $1\leq s<\infty$.\par\noindent
We set, for $\eps_0>0$,
$$\Omega_{\eps_0}=\{x\in S^n:\ |v_\eps(x)-\bar c_n|<\delta\quad\hbox{for}\quad \eps\leq\eps_0\}.$$
Since $v_\eps\to \bar c_n$ in measure, one has
$$|S^n\,\backslash\, \Omega_{\eps_0}|\longrightarrow 0,\quad\hbox{as}\quad \eps_0\to 0.$$
Let us write problem \ellegizeroprimo\ in the form
$$\laplconf^2(v_\eps-\bar c_n)=q_\eps(v_\eps-\bar c_n)+w_\eps$$
where
$$\eqalign{&q_\eps=b_\eps\chi_{S^n\,\backslash \Omega_{\eps_0}}\quad\hbox{and}\cr
&w_\eps=b_\eps(v_\eps-\bar c_n)\chi_{\Omega_{\eps_0}}+\bar c_n\left(b_\eps-{c_n^2\over 4}\right).\cr}$$
From Lemma 1.3, there exists $c_1(n)>0$ such that
$$\|q_\eps\|_{L^\infty}\leq c_1(n),\quad \|w_\eps\|_{L^\infty}\leq c_1(n),\quad \|q_\eps\|_{L^{n\over 4}}\leq c_1(n)|S^n\,\backslash \Omega_{\eps_0}|^{4\over n}.$$
Using now Proposition 1.4 with $p=1$, $q=q_\eps$, $g=w_\eps$ and $s>{n\over4}$, we obtain that there exists a positive constant $c_2(n)$ and $\eps_1>0$ small enough such that for any $\eps\leq\eps_1$
$$\|v_\eps-\bar c_n\|_\infty\leq c_2(n)\|w_\eps\|_{L^s}.\eqno\stimauno$$
Since $b_\eps\to {c_n^2\over 4}$ in $L^s(S^n)$, by definition of $\Omega_{\eps_0}$, it is clear that there exists $\eps_2$ such that for any $\eps\leq \eps_2$
$$\|w_\eps\|_{L^s}\leq c_3\delta\eqno\stimadue$$
for some positive constant $c_3$ independent of $\delta$ and $\eps$. Since $\delta$ can be chosen arbitrarily small, from \stimauno\ and \stimadue\ we have that, for $\eps$ small, $v_\eps$ is closed in the uniform norm to a positive number; in particular $v_\eps$ (and so $u_\eps$) is positive.\fine
\medskip\noindent
{\bf Theorem 1.5.}\quad\sl If $K$ satisfy (K1-2). Then for $\eps$ small \tre\ admits at least one positive solution.\rm
\medskip\noindent
{\bf Proof.}\quad From Lemma 1.2 $f_\eps$ has a critical point $u_\eps\in E$ solving \quattro. The preceding arguments imply that $u_\eps$ is positive and hence it is a solution of \tre, too.\fine
\medskip\noindent
{\bf 1.3. Other results.}\quad In this subsection we discuss some results which are the counterpart of those of [\spagnoli, Sections 4,5,6]. For the sake of brevity, we will state here only some specific results. The proofs are easily obtained, as before, extending all the arguments of [\spagnoli] to our setting. Of course, the positivity of the solutions has to be proved by means of the Van Der Vorst techniques, as carried out in the preceding subsection.
\medskip\noindent
{\bf Theorem 1.6.} \quad \sl Let $K$ satisfy
$$\left\{\eqalign{&K\in L^{\infty}(\erre^n)\cap C^1(\erre^n),\cr
                  &K(y)=K(r),\quad r=|y|,\cr
                  &r^{-\alpha}K(r)\in L^1([0,\infty),r^{n-1}\,dr)\quad\hbox{for some}\quad \alpha<n\cr}\right.\leqno{(K3)}$$
 and $K(0)=0$, $K\not\equiv0$. Then for $|\eps|$ small, \tre\ has a positive radial solution $u_\eps\in{\cal D}^{2,2}$.
\medskip\noindent \rm
{\bf Theorem 1.7.} \quad \sl \smallskip\noindent
\item{$\bullet$} Let $h$ satisfy 
$$\hbox{(h1.a)}\quad h\in C^2(\erre^n)\cap  L^{\infty}(\erre^n)\hskip10truecm$$
and $h\not\equiv 0$. Then for $|\eps|$ small
$$\Delta^2 u=\eps h(y) u+u^p,\quad y\in\erre^n,\quad n>8\eqno\tredue$$
has a positive solution $u_\eps\in{\cal D}^{2,2}(\erre^n)$. Furthermore, if $\exists\, \xi_1,\xi_2\in\erre^n$ such that $h(\xi_1)>0$, $h(\xi_2)<0$, then for $|\eps|$ small \tredue\ has at least two distinct positive solutions. \smallskip\noindent
\item{$\bullet$} Let $h\in L^1(\erre^n)\cap L^{\infty}(\erre^n)$ with $\int_{\erre^n}h\not=0$. Then for $|\eps|$ small
$$\Delta^2u=\eps h(y) u+u^p,\quad y\in\erre^n,\quad n\geq 5$$
has a positive solution $u_\eps\in{\cal D}^{2,2}(\erre^n)$.
\medskip\noindent\rm
{\bf Theorem 1.8.} \quad \sl Let us assume $n\geq5$, $1<q<p$, and $h\in L^1(\erre^n)\cap L^{\infty}(\erre^n)$, $h\not\equiv 0$. Then for $|\eps|$ small enough problem 
$$\Delta^2u=\eps h(y)u^q+u^p,\quad y\in\erre^n,\quad 1<q<p={n+4\over n-4}\eqno\quattrouno$$
has a positive solution in ${\cal D}^{2,2}(\erre^n)$.\medskip\noindent
\rm

\vskip2truecm
\nuovoparagrafo 2

\centerline{\bf 2.\quad The non-perturbative case}
\bigskip\noindent
In order to treat the non-perturbative case we follow the procedure that Y. Y. Li has used for the scalar curvature problem in [\liuno]. The present section is divided in three parts; in the first we make a blow up analysis for $(P)$ and prove that only simple isolated blow up points can occur. Then we show that under some flatness assumption there are not blow up points and finally we make a homotopy in order to pass from the perturbative case to the non perturbative one.
\medskip\noindent
{\bf 2.1. Blow up points have to be simple isolated.}\quad  We are interested in the family of problems
$$\pan v={n-4\over 2}\widehat Q_i(x)v^p,\qquad v>0,\quad\hbox{on}\ S^n.\eqno\seitre_i$$
If $v_i$ is a solution of $\seitre_i$ then $\ui(y)=\varphi(y)v_i(\pi^{-1}(y))$ is a solution of
$$\Delta^2 u={n-4\over 2}Q_i(y)u^p,\quad u>0,\quad\hbox{in}\ \erre^n,\eqno\seiquattro_i$$
where $\varphi$ is given by \dueprimo\ and $Q_i=\widehat Q_i\circ\pi^{-1}$. The question is what happens to $\ui$ when $i\to\infty$. For the reader's convenience here we report the definition of blow up point that one can find in [\dmadue]. Let $B_2=\{y\in\erre^n:\ |y|<2\}$ and $Q_i$ be a sequence of $C^1(B_2)$ functions such that
$${1\over A_1}\leq Q_i(y)\leq A_1\qquad\forall\, i,\quad\forall\,y\in B_2,\quad\hbox{for some}\ A_1>0\eqno\seiuno$$
and consider the equations
$$\Delta^2 u={n-4\over 2}Q_i(y)u^p,\quad u>0,\quad -\Delta u>0,\quad y\in B_2.\eqno\seidue_i$$
\medskip\noindent
{\bf Definition 2.1.}\quad \sl Let $\{u_i\}_i$ be a sequence of solutions of $\seidue_i$. A point $\bar y\in B_2$ is called a blow up point of $\{u_i\}_i$ if there exists a sequence $\{y_i\}_i$ converging to $\bar y$ such that $u_i(y_i)\to\infty$. 
\medskip\noindent
\rm In the sequel, if $\bar y$ is a blow up point for $\{\ui\}_i$, writing $\yi\to \bar y$ we mean that $\yi$ are local maxima of $\ui$ and $\ui(\yi)\to+\infty$ as $i\to+\infty$.
\medskip\noindent
{\bf Definition 2.2.}\quad \sl A blow up point $\bar y$ is said to be an isolated blow up point of $\{u_i\}_i$ if there exist $\bar r\in(0, {\rm dist}\, (\bar y,\partial B_2))$, $\bar c>0$ and a sequence $\{y_i\}$ converging to $\bar y$ such that $y_i$ is a local maximum of $u_i$, $u_i(y_i)\to\infty$ and
$$u_i(y)\leq\bar c|y-y_i|^{-{4\over p-1}}\quad\forall\, y\in B_{\bar r}(y_i).$$
\medskip\noindent
{\bf Definition 2.3.}\quad An isolated blow up point $\bar y\in B_2$ of $\{u_i\}_i$ is said to be an isolated simple blow up point of $\{u_i\}_i$ if there exists some positive $\rho$ independent of $i$ such that
$$\bar w_i(r)={r^{4\over p-1}\over|\partial B_r|}\int_{\partial B_r} u_i\qquad r>0$$
has precisely one critical point in $(0,\rho)$ for large $i$.
\medskip\noindent\rm
Suppose that there exist 
$$\left\{\eqalign{&\beta\geq2,\quad\beta\in(n-4,n),\cr
                  &A_1,L_1,L_2>0,\cr}\right.$$
such that
$$\leqalignno{&Q_i\in C_{{\rm loc}}^{[\beta]-1,1}(B_2)\quad \hbox{and}{1\over A_1}\leq Q_i\leq A_1\quad \hbox{in}\ B_2;&(Q1)\cr
              &\|\nabla Q_i\|_{C^0(\Omega_i)}\leq L_1\quad\hbox{and}&(Q2)\cr
              & |\nabla^sQ_i(y)|\leq L_2|\nabla Q_i(y)|^{\beta-s\over\beta-1}\quad\forall\, 2\leq s\leq [\beta],\quad\forall\, y\in\Omega_i:\ \nabla Q_i(y)\not=0,&\cr}$$
where $\Omega_i$ is a neighbourhood of the set of critical points of $Q_i$.\par\noindent
In the sequel $c_1,\ c_2, \dots$ will denote positive constants which may vary from formula to formula and which may depend only on $\rho,A_1,L_1,L_2,n$. The following two lemmas are proved in the Appendix following the scheme of Lemma 2.6 and 2.8 in [\liuno] and using some results of [\dmadue].
\medskip\noindent
{\bf Lemma 2.4.}\quad\sl Assume that (Q1-2) hold in $B_2$. Let $u_i$ satisfy $\seidue_i$ and suppose that $\bar y$ $(\yi\to\bar y)$ is an isolated simple blow up point. We have that
$$|\nabla Q_i(\yi)|\leq c_1\ui^{-2}(\yi)+c_2\ui^{-{2\over n-4}\left({{\beta-[\beta]\over\beta-1}+[\beta]-1}\right)}(\yi).$$
\medskip\noindent
{\bf Lemma 2.5.}\quad\sl Under the same assumptions as in Lemma 2.4 and for all $0<\sigma<1$ we get:\par for $\beta=2$
$$\left|\int_{B_\sigma(\yi)}(y-\yi)\cdot\nabla Q_i\ui^{\pii+1}\right|
\leq c\ui^{-{4\over n-4}}(\yi)+o\left(\nabla Q_i(\yi)\ui^{-{2\over n-4}}(\yi)\right)$$
\par and for $\beta>2$
$$\left|\int_{B_\sigma(\yi)}(y-\yi)\cdot\nabla Q_i\ui^{\pii+1}\right|\leq c_1|\nabla Q_i(\yi)|\ui^{-{2\over n-4}}(\yi)+c_2\ui^{- {2\over n-4}\left([\beta]+{\beta-[\beta]\over\beta-1}\right)}(\yi).$$
\medskip\noindent \rm Thanks to the previous lemmas, we obtain a corollary, which is the analogous of Co\-rol\-lary~2.1 of [\liuno].
\medskip\noindent
{\bf Corollary 2.6.}\quad\sl  Assume (Q1-2). Let $u_i$ satisfy $\seidue_i$ and suppose that $y_i\to0$ is an isolated simple blow up point. Then for any $0<\sigma<1$ we have that
$$\lim_{i\to\infty}\ui^2(\yi)\int_{B_\sigma(\yi)}(y-y_i)\cdot\nabla Q_i\ui^{\pii+1}=0.$$
\medskip\noindent
{\bf Proof.}\quad\rm Owing to Lemma 2.4 and Lemma 2.5 we have that
$$\eqalign{&\left|\ui^2(\yi)\int_{B_\sigma(\yi)}(y-y_i)\cdot\nabla Q_i\ui^{\pii+1}\right|\cr
&\qquad\leq c_1\ui^{2-{2\over n-4}}(\yi)\left[c_2\ui^{-2}(\yi)+c_3\ui^{-{2\over n-4}\left([\beta]-1+{\beta-[\beta]\over\beta-1}\right)}(\yi)\right]\cr
&\qquad=c_4\ui^{-{2\over n-4}}(\yi)+c_5\ui^{2-{2\over n-4}\left([\beta]+{\beta-[\beta]\over\beta-1}\right)}(\yi).\cr}$$
Remark that 
$[\beta]+{\beta-[\beta]\over\beta-1}>n-4$
to obtain the statement.\fine
\medskip\noindent
Corollary 2.6 can be used to prove that, in our hypotheses, blow up points for equation $\seidue_i$  have to be {\sl isolated simple} blow up points. In particular, note that we can prove the same results of [\dmadue, Proposition 2.19 and 3.2] (where (Q2) was not required) for any $n\geq5$.
\medskip\noindent
{\bf Theorem 2.7.}\quad \sl In the same assumptions of Lemma 2.4, blow up points for equation $\seidue_i$ are isolated simple and finite.
\medskip\noindent\rm
{\bf Sketch of the proof.}\par {\bf Step 1. Any isolated blow up point has to be isolated simple.}\quad We can proceed as in [\dmadue, Proposition 2.19]. To prove that any isolated blow up point $\bar y$ for $\{\ui\}_i$ has to be simple, in [\dmadue, Proposition 2.19] it is shown that the function $\bar w_i$ given by Definition 2.3 has precisely one critical point in $(0,R_i\ui^{-(p-1)/4})$ (for a suitable $R_i\to+\infty$). By contradiction suppose that $\bar y$ is not simple and let $\mu_i\geq R_i\ui^{-{p-1\over4}}$ be the second critical point of $\bar w_i$. Assume $\yi=0$ and set $\xi(y)=\mu_i^{4/(p-1)}\ui(\mu_iy)$,\quad $|y|\leq 1/\mu_i$. $\xi_i$ satisfies
$$\left\{\eqalign{&\Delta^2\xi_i(y)=Q_i(\mu_i y)\xi_i(y)^p,\quad\qquad |y|<{1\over\mu_i}\cr
                 &|y|^{4\over p-1}\xi_i(y)\leq {\rm const},\ \quad\qquad\quad\quad|y|<{1\over\mu_i}\cr
&\lim_i\xi_i(0)=+\infty.\cr}\right.\eqno\semplice$$
In [\dmadue] it is proved that $0$ is an isolated simple blow up point for equation \semplice\ and that
$$\xi_i(0)\xi(y)\longrightarrow h(y)=a|y|^{4-n}+b(y)\qquad\hbox{in}\quad C^4_{\rm loc}(\erre^n\,\backslash\,{0})$$
for some biharmonic function $b$ and for some constant $a>0$. Then it is shown that setting
$$\eqalignno{B(r,y,h,\nabla h,\nabla^2h,\nabla^3 h)&=-{n-2\over 2}\Delta h{\partial h\over\partial \nu}-{r\over 2}|\Delta h|^2+{n-4\over 2}h{\partial\over\partial\nu}(\Delta h)&\cr
  &+\dual{y,\nabla h}{\partial\over\partial\nu}(\Delta h)-\Delta h\sum_iy_i{\partial\over\partial\nu}h_i.&\defdibi\cr}$$
one has
$$0>\int_{\partial B_{\sigma}}B(\sigma,y,h,\nabla h,\nabla^2h,\nabla^3h)\geq \lim_i \xi_i(0)^2{n-4\over 2(p+1)}\sum_j\int_{B_\sigma}y_j{\partial Q_i(\mu_i\cdot)\over\partial y_j}\xi_i^{p+1}.\eqno\defdibiprimo$$
Only at this point of the proof of [\dmadue] the fact that $n=5,6$ is used to conclude that the last term of the previous inequality is $0$ thus finding a contradiction. However, in our case it works for every $n\geq5$ thanks to Corollary 2.6.
\smallskip{\bf Step 2. Blow up points have to be isolated.}\quad We observe that, using Corollary 2.6 again, we are able to extend [\dmadue, Proposition 3.2] which states that for any $p,q\in S^n$ blow up points of $\seidue_i$ one has
$$|p-q|\geq\delta^*\eqno\finito$$
under flatness condition. Hence the number of blow up points is finite. \finito\ and [\dmadue, Proposition 3.1, statement (3)] show that the blow up points are indeed isolated. Therefore we get the thesis of the theorem.\fine
\medskip\noindent
{\bf 2.2. Blow up analysis under flatness assumption.} The following theorem ensures that under our hypotheses there is at most one isolated simple blow up point. We prove it following the scheme of the proof of [\liuno, Theorem 4.2].
\medskip\noindent
{\bf Theorem 2.8.}\quad\sl Suppose that $\{\widehat Q_i\}_{i\in\enne}\subset C^1(S^n)$ with uniform $C^1$ modulo of continuity and
$${1\over A_2}\leq \widehat Q_i(q)\leq A_2\quad\forall\, q\in S^n.$$
Suppose also that there exists $\beta\in(n-4,n),$ $\beta\geq 2$ such that
$$\left\{\eqalign{&\widehat Q_i\in C^{[\beta]-1,1}(\Omega_{d,i}),\qquad\|\nabla_{g_0} \widehat Q_i\|_{C^0(\Omega_{d,i})}\leq L_1,\cr  
                  &|\nabla^s_{g_0}\widehat Q_i(y)|\leq L_2|\nabla_{g_0}\widehat Q_i(y)|^{\beta-s\over\beta-1}\quad\forall\, 2\leq s\leq[\beta],\ y\in \Omega_{d,i},\ \nabla_{g_o}\widehat Q_i(y)\not= 0,\cr}\right.\leqno(*)_\beta$$
where 
$$\Omega_{d,i}=\{q\in S^n:\ |\nabla_{g_0}\widehat Q_i(q)|<d\}.$$
Let $v_i$ be solutions of $\seitre_i$. Then, after passing to a subsequence, either $\{v_i\}_i$ stays bounded in $L^\infty(\Omega)$ (and hence in $C^{4,\alpha}(S^n),\ 0<\alpha<1$) or $\{v_i\}$ has precisely one isolated simple blow up point.
\medskip\noindent
{\bf Proof.}\quad \rm Owing to Theorem 2.7, any blow up point has to be simple isolated. Suppose by contradiction that $\{v_i\}_i$ has two distinct simple isolated blow up points $q^1\not= q^2$ and let $q^1_i$, $q^2_i$ be local maxima of $v_i$ with
$$q^1_i\to q^1\quad\hbox{and}\quad q^2_i\to q^2.$$
We may assume that $q^2\not=- q^1$ without loss of generality. We make a stereographic projection with $q^1_i$ as the south pole. Note that if $q$ is a blow up point for $\{v_i\}_i$ then $\pi(q)$ is a blow up point for $\{\ui\}_i=\{v_i\circ\pi^{-1}\}_i$. In the stereographic coordinates $q^1_i=0$ and we still use $q^2$, $q^2_i$ to denote the stereographic coordinates of $q^2$, $q^2_i\in S^n$. Then we get equation $\seiquattro_i$. It is easy to see that $\{Q_i\}_i$ satisfies $(Q2)$ for some constants $L_1$ and $L_2$ in some open set of $\erre^n$ containing $0$ and $q^2$. 
Then an application of Lemma 2.4 yields $\nabla Q_i(q^{1,2}_i)\to 0$. \par
\finito\  implies that the number of blow up points is bounded by some constant independent of $i$. Therefore we can argue as in the sketch of the proof of Theorem 2.7 and see that there exist some finite set ${\cal S}\subseteq\erre^n$, $0,q^2\in {\cal S}$, some constants $a, A>0$ and some function $h(y)\in  C^4(\erre^n\,\backslash\,{\cal S})$ such that
$$\eqalignno{&\lim_iu_i(0)u_i(y)=h(y)\quad\hbox{in}\quad C^4_{{\rm loc}}(\erre^n\,\backslash\, {\cal S}),&\punto\cr
&h(y)=a|y|^{4-n}+A+o(|y|),\quad\hbox{for}\quad y\sim 0.&\cr}$$
Corollary 2.6 implies that for any $0<\sigma<1$
$$\lim_{i}\ui(0)^2\int_{B_\sigma}y\cdot\nabla Q_i\ui^{\pii+1}=0.\eqno\contradue$$ 
Corollary 4.2 (part (ii)) in [\dmadue] yields
$$\lim_{r\to 0}\int_{\partial B_r}B(r,y,h,\nabla h,\nabla^2h,\nabla^3 h)=-(n-4)^2(n-2)|S^{n-1}|<0\eqno\contrauno$$
where $B$ is defined in \defdibi. On the other side, using the Pohozaev type arguments one finds, see [\dmadue, Proposition 4.1],
$$\eqalignno{\int_{\partial B_\sigma}B(\sigma,y,\ui,\nabla\ui,\nabla^2\ui,\nabla^3\ui)=&{n-4\over 2(p+1)}\int_{B_\sigma}y\cdot\nabla Q_i(y)\ui^{p+1}&\cr
&-\sigma{n-4\over 2(p+1)}\int_{\partial B_\sigma}Q_i\ui^{p+1}.&\poho\cr}$$
Recall now that from [\dmadue, Proposition 2.9] 
$$C^{-1}{|y-\yi|^{4-n}\over\ui(\yi)}\leq\ui(y)\leq C{|y-\yi|^{4-n}\over\ui(\yi)}
\quad\hbox{for}\quad R_i\ui^{-{p-1\over4}}(\yi)\leq|y-\yi|\leq 1\eqno\propduenove$$
for a positive constant $C$ and for any $R_i\to\infty$ (up to a subsequence).
Multiplying \poho\ by $\ui(0)^2$, using \punto\ and noting that from \propduenove
$$\int_{\partial B_\sigma}Q_i\ui^{p+1}\leq c \ui^{-p-1}(\yi)\mathop{\longrightarrow}\limits_{i\to\infty} 0$$
we get
$$\eqalign{\int_{\partial B_\sigma}B(\sigma,y,h,\nabla h,\nabla^2h,\nabla^3 h)&=\lim_i \ui(0)^2\int_{\partial B_\sigma}B(\sigma,y,\ui,\nabla \ui,\nabla^2\ui,\nabla^3 \ui)\cr
&= \lim_i \ui(0)^2{n-4\over 2(\pii+1)}\sum_j\int_{B_\sigma}y_j{\partial Q_i\over \partial y_j}\ui^{\pii+1}.\cr}$$
From \contrauno\ we know that, for $\sigma$ small enough, the first integral in the above formula is strictly less than 0, whereas \contradue\ implies that the last term above is 0, which is a contradiction.\fine
\medskip\noindent
We will now prove that if a point $q_0\in S^n$ satisfies the following condition
$$\left\{\eqalign{&\exists\, \eps_0>0:\ \widehat Q_i\in C^{[\beta]-1,1}(B_{\eps_0}(q_0))\quad\hbox{and}\cr
                  &\widehat Q_i(q_0)\geq{1\over A_2}\quad\hbox {for some positive constant }A_2\cr
&\widehat Q_i(y)=\widehat Q_i(q_0)+\sum_{j=1}^na_j|y_j|^\beta+{\cal R}_i(y),\quad |y|\leq \eps_0,\cr
&\hbox{where }a_j=a_j(q_0)\not=0,\quad\sum_{j=1}^na_j\not=0,\cr
&y\ \hbox{is some geodesic normal coordinate system centered at}\ q_0,\cr
& {\cal R}_i(y)\in C^{[\beta]-1,1}\ \hbox{near}\ 0\ \hbox{and}\cr
&\lim_{|y|\to 0}\sum_{0\leq s\leq[\beta]}|\nabla^s{\cal R}_i(y)||y|^{-\beta+s}=0\quad\hbox{uniformly for}\ i,\cr}\right.\leqno(**)_\beta$$
then the form of $\widehat Q_i$ near $q_0$ prevents blow up from happening there. The following theorem is the analogous of [\liuno, Corollary 4.2].
\medskip\noindent 
{\bf Theorem 2.9.}\quad \sl Suppose $\{\widehat Q_i\}_i\subseteq C^1(S^n)$ with uniform $C^1$ modulo of continuity and satisfies for some $q_0\in S^n$ condition $(**)_\beta$. Let $v_i$ be positive solutions of $(P)$ with $\widehat Q=\widehat Q_i$. If $q_0$ is an isolated simple blow up point of $v_i$, then $v_i$ has to have at least another blow up point.
\medskip\noindent
\rm\medskip\noindent
{\bf Remark 2.10.}\quad Note that $(**)_\beta$ implies that $\{Q_i\}_i=\{\widehat Q_i\circ \pi^{-1}\}_i$ satisfies (Q2). See [\liuno, Remark 0.3].
\rm\medskip\noindent
{\bf  Proof of Theorem 2.9.}\quad Let us argue by contradiction. Suppose that $v_i$ has no other blow up points. Let us take $q_0$ as the south pole and make stereographic projection on $\erre^n$, so that our equation becomes
$$\left\{\eqalign{\Delta^2 \ui={n-4\over2}Q_i(y)\ui(y)^{n+4\over n-4}\quad&\hbox{in}\ \erre^n,\cr
\ui>0 \quad&\hbox{in}\ \erre^n,\cr}\right.\eqno\seidiciannove$$
where $\ui(y)=\varphi(y) v_i(\pi^{-1}(y)).$ Let $\yi$ be the local maximum of $\ui$ and such that $\yi\to 0$. It follows from Lemma 2.4 that
$$|\nabla Q_i(\yi)|=O\left(\ui^{-2}(\yi)+\ui^{ -{2\over n-4}\left([\beta]-1+{\beta-[\beta]\over\beta-1}\right)}(\yi)\right).$$
Multiply equation \seidiciannove\  by $\nabla\ui$; it is easy to check by integration by parts that $$\int_{\erre^n}\Delta^2\ui\nabla\ui=0.$$ Then we obtain
$$\int_{\erre^n}Q_i\ui^{n+4\over n-4}\nabla\ui=0.$$
Another integration by parts yields
$$0=\int_{\erre^n}Q_i\ui^{n+4\over n-4}\nabla\ui =-{n-4\over 2n}\int_{\erre^n}\ui^{2n\over n-4}\nabla Q_i.\eqno\seiventi$$
Since we have assumed that $v_i$ has no other blow up point but $q_0$, the Harnack inequality and \propduenove\ yield that
$$v_i(q)\sim(\max_{S^n}v_i)^{-1}\quad\hbox{for}\ q\ \hbox{away from}\ q_0,$$
and this implies that for $|y|\geq\eps>0$ $\ui(y)\leq c(\eps)|y|^{4-n}\ui^{-1}(\yi)$.\par\noindent
Therefore for $\eps>0$ small from \seiventi\ we obtain
$$\eqalignno{&\left|\int_{B_\eps}\nabla Q_i(y+\yi)\ui^{2n\over n-4}(y+\yi)\, dx\right|\leq\left|\int_{\erre^n}\nabla Q_i(y+\yi)\ui^{2n\over n-4}(y+\yi)\, dx\right|&\cr
&\qquad+\left|\int_{\erre^n\,\backslash\,B_\eps}\nabla Q_i(y+\yi)\ui^{2n\over n-4}(y+\yi)\, dx\right|\leq C(\eps)\ui^{-{2n\over n-4}}(\yi).&\seiventuno\cr}$$
Let us set ${\cal T}^{(\beta)}(y)=\sum_j^n a_j|y_j|^\beta.$
Using $(**)_\beta$ \seiventuno\ becomes
$$\left|\int_{B_\eps}[\nabla {\cal T}^{(\beta)}(y+\yi)+\nabla {\cal R}_i(y+\yi)]\ui^{2n\over n-4}(y+\yi)\, dx\right|\leq C(\eps)\ui^{-{2n\over n-4}}(\yi).\eqno\seiventunoprimo$$
Since
$${|\nabla {\cal R}_i(y)|\over|\nabla {\cal T}^{(\beta)}(y)|}\leq c{|\nabla {\cal R}_i(y)|\over|y|^{\beta-1}}\mathop{\longrightarrow}\limits_{y\to 0}0\quad\hbox{uniformly in}\ i,$$
from \seiventunoprimo \ we deduce
$$\left|\int_{B_\eps}\left\{\nabla {\cal T}^{(\beta)}(y+\yi)+o_\eps(1)|\nabla {\cal T}^{(\beta)}(y+\yi)|\right\}\ui^{2n\over n-4}(y+\yi)\right|\leq C(\eps)\ui^{-{2n\over n-4}}(\yi)$$
where $o_\eps(1)$ stays for something going to $0$ as $\eps\to 0$ uniformly in $i$. Multiplying the above by $\ui^{{2\over n-4}(\beta-1)}(\yi)$, noting that $\nabla{\cal T}^{(\beta)}(\lambda y)=\lambda^{\beta-1}\nabla{\cal T}^{(\beta)}(y)$ for any $\lambda>0$, and setting
$$\eqalign{&\xi_i:=\ui^{2\over n-4}(\yi)\yi\cr
           &A_i=\int_{B_\eps}\bigg\{\nabla {\cal T}^{(\beta)}\Big(\ui^{2\over n-4}(\yi)y+\xi_i\Big)+o_\eps(1)\left|\nabla {\cal T}^{(\beta)}\Big(\ui^{2\over n-4}(\yi)y+\xi_i\Big)\right|\bigg\}\ui^{2n\over n-4}(y+\yi),\cr}$$
we get
$$|A_i|\leq C(\eps)\ui^{{2\over n-4}(\beta-1-n)}(\yi).\eqno\ai$$
We now want to prove that $\{\xi_i\}_i$ is bounded. By contradiction, suppose that $\xi_i\mathop{\to}\limits_{i\to\infty}\infty$ along a subsequence. Let us set
$$\eqalign{I_i=\int_{|y|\leq R_i\ui^{-{2\over n-4}}(\yi)}\bigg\{&\nabla {\cal T}^{(\beta)}\Big(\ui^{2\over n-4}(\yi)y+\xi_i\Big)\cr
&+o_\eps(1)\left|\nabla {\cal T}^{(\beta)}\Big(\ui^{2\over n-4}(\yi)y+\xi_i\Big)\right|\bigg\}\ui^{2n\over n-4}(y+\yi)\cr}$$
and
$$\eqalign{J_i=\int_{R_i\ui^{-{2\over n-4}}(\yi)\leq |y|\leq \eps}\Big\{&\nabla {\cal T}^{(\beta)}\left(\ui^{2\over n-4}(\yi)y+\xi_i\right)\cr
&+o_\eps(1)\left|\nabla {\cal T}^{(\beta)}\left(\ui^{2\over n-4}(\yi)y+\xi_i\right)\right|\Big\}\ui^{2n\over n-4}(y+\yi).\cr}$$
Making the change of variable $z=\ui^{2\over n-4}(\yi)y$ and using (4.7) of Appendix we get
$$|I_i|\sim c\Bigg|\int_{|z|\leq R_i}(z+\xi_i)|z+\xi_i|^{\beta-2}\left(\ui^{-1}(\yi)\ui\Big(\ui^{-{2\over n-4}}(\yi)z+\yi\Big)\right)^{2n\over n-4}\Bigg|\sim|\xi_i|^{\beta-1}\eqno\ii$$
and 
$$|J_i|\leq c\left|\int_{R_i\ui^{-{2\over n-4}}(\yi)\leq |y|\leq \eps}\Big\{|\ui^{2\over n-4}(\yi)y|^{\beta-1}+|\xi_i|^{\beta-1}\Big\}\ui(y+\yi)^{2n\over n-4}\right|.$$
Here and in the sequel we need the following estimate which one can find in [\dmadue, Lemma 2.17]: in the assumptions of Lemma 2.4 and for any sequence $R_i\to+\infty$ one has that (up to a subsequence)
$$\left\{\eqalign{&\int_{B_{r_i}}|y-\yi|^s\ui^{p+1}=O\left(\ui^{-{2s\over n-4}}(\yi)\right)\cr
&\int_{B_1\,\backslash\,B_{r_i}}|y-\yi|^s\ui^{p+1}\leq o\left(\ui^{-{2s\over n-4}}(\yi)\right)\cr}\right.\qquad -n<s<n\eqno\stima_s$$
for some positive constant $c$, where $r_i=R_i\ui^{-{p-1\over 4}}(\yi)$. The second of $\stima_s$ with $s=\beta-1$ gives 
$$|J_i|\leq o\left(\ui^{{2(\beta-1)\over n-4}-{2\over n-4}(\beta-1)}(\yi)\right)+|\xi_i|^{\beta-1}o(1) \leq o(1)|\xi_i|^{\beta-1}.\eqno\ji$$
Putting together \ai, \ii, and \ji\ we find
$$ C(\eps)\ui^{{2\over n-4}(\beta-1-n)}(\yi)\geq|A_i|\geq|I_i|-|J_i|\geq|\xi_i|^{\beta-1}(1-o(1))\sim|\xi_i|^{\beta-1}.$$
Then
$$|\xi_i|\leq c(\eps)\ui^{{2\over n-4}\left(1-{n\over\beta-1}\right)}(\yi)\leq c$$
for some positive $c$, which is a contradiction. Then 
$$|\yi|=\ui^{-{2\over n-4}}(\yi)|\xi_i|=O\left(\ui^{-{2\over n-4}}(\yi)\right).\eqno\seiventidue$$
Writing \poho\ in $\erre^n$ for equation \seidiciannove\ we have that
$${n-4\over 2(p+1)}\int_{\erre^n}y\cdot\nabla Q_i(y)\ui^{p+1}=0.$$
Using \seiventi, we get
$$\int_{\erre^n}y\cdot \nabla Q_i(y+\yi)\ui^{2n\over n-4}(y+\yi)=0$$
so that for any $\eps>0$ small we have 
$$\left |\int_{B_\eps}y\cdot \nabla Q_i(y+\yi)\ui^{2n\over n-4}(y+\yi)\right|\leq c(\eps)\ui^{-{2n\over n-4}}(\yi).$$
Therefore $\stima_s$ with $s=\beta$ and $s=1$ and $(**)_\beta$ yield
$$\eqalign{&\int_{B_\eps} y\cdot \nabla {\cal T}^{(\beta)}(y+\yi)\ui^{2n\over n-4}(y+\yi)\cr
&\qquad \leq c(\eps)\ui^{-{2n\over n-4}}(\yi)+o_\eps(1)\int_{B_\eps}|y||y+\yi|^{\beta-1}\ui^{2n\over n-4}(y+\yi)\cr
 &\qquad \leq c(\eps)\ui^{-{2n\over n-4}}(\yi)+o_\eps(1)\int_{B_\eps}\{|y|^\beta+|y||\yi|^{\beta-1}\}\ui^{2n\over n-4}(y+\yi)\cr
&\qquad \leq c(\eps)\ui^{-{2n\over n-4}(\yi)}+o_\eps(1)\ui^{-{2\beta\over n-4}}(\yi)+o_\eps(1)|\yi|^{\beta-1}\ui^{-{2\over n-4}}(\yi).\cr}$$
Multiplying the above by $\ui^{2\beta\over n-4}(\yi)$ and using \seiventidue\ we have that
$$\lim_{i\to\infty}\ui^{2\beta\over n-4}(\yi)\left|\int_{B_\eps}y\cdot \nabla {\cal T}^{(\beta)}(y+\yi)\ui^{2n\over n-4}(y+\yi)\, dy\right|\leq o_\eps(1).\eqno\seiventitre$$
Moreover $\stima_s$ with $s=\beta$ and $s=1$ and \seiventidue\ yield
$$\eqalignno{&\ui^{2\beta\over n-4}(\yi)\left|\int_{R_i\ui^{-{2\over n-4}}(\yi)\leq |y|\leq\eps}y\cdot \nabla {\cal T}^{(\beta)}(y+\yi)\ui^{2n\over n-4}(y+\yi)\, dy\right|&\cr
&  \leq c\ui^{2\beta\over n-4}(\yi)\Big\{\int_{R_i\ui^{-{2\over n-4}}(\yi)\leq|y|\leq\eps}(|y|^\beta+|y||\yi|^{\beta-1})\ui^{2n\over n-4}(y+\yi)\, dy\Big\}\mathop{\longrightarrow}\limits_{i\to\infty}0.\quad\quad&\seiventiquattro\cr}$$
Set $\xi:=\lim_{i\to\infty}\ui^{2\over n-4}(\yi)\yi$ (see \seiventidue). Making the change of variable $y=\ui^{-{2\over n-4}}(\yi)z$ and using \accorpa\seiventitre\seiventiquattro\ we get
$$\left|\int_{|z|\leq R_i}z\cdot\nabla{\cal T}^{(\beta)}(z+\xi_i)\left[\ui^{-1}(\yi)\ui\left(\ui^{-{2\over n-4}}z+\yi\right)\right]^{2n\over n-4}\,dz\right|=o_\eps(1)+o(1).$$
Thanks to (4.7) of Appendix we have that
$$\ui^{-1}(\yi)\ui\left(\ui^{-{2\over n-4}}z+\yi\right)\mathop{\longrightarrow}\limits_{i\to\infty}(1+k|z|^2)^{-n}\quad\hbox{uniformly}$$
where $k^2=\lim_{i\to\infty}[2n(n-2)(n+2)]^{-1}Q_i(\yi)$.
Therefore letting $i\to\infty$ we get
$$\left|\int_{\erre^n}z\cdot\nabla {\cal T}^{(\beta)}(z+\xi)(1+k|z|^2)^{-n}\, dz\right|=o_\eps(1).$$
Hence we have that
$$\int_{\erre^n}z\cdot \nabla {\cal T}^{(\beta)}(z+\xi)(1+k|z|^2)^{-n}\, dz=0.\eqno\seiventicinque$$
Using \seiventi\ and arguing as before we have
$$\int_{\erre^n} \nabla {\cal T}^{(\beta)}(z+\xi)(1+k|z|^2)^{-n}\, dz=0.\eqno\seiventisei$$
Using \accorpa\seiventicinque\seiventisei\ and noting that ${\cal T}^{(\beta)}(x)={1\over\beta}x\cdot \nabla {\cal T}^{(\beta)}(x)$, we deduce
$$\int_{\erre^n}{\cal T}^{(\beta)}(z+\xi)(1+k|z|^2)^{-n}\,dz=\beta^{-1}\int_{\erre^n}(z+\xi)\cdot\nabla {\cal T}^{(\beta)}(z+\xi)(1+k|z|^2)^{-n}\, dz=0.\eqno\seiventisette$$
For any $\teta\in\erre^n$ we have that
$$\int_{\erre^n}{\partial\over\partial y_j}{\cal T}^{(\beta)}(y+\teta)(1+k|y|^2)^{-n}\,dy=\beta a_j\int_{\erre^n}|y_j+\teta_j|^{\beta-2}(y_j+\teta_j)(1+k|y|^2)^{-n}\,dy.$$
It is easy to prove that the last integral is $0$ if and only if $\teta_j=0$, so that
$$\int_{\erre^n}\nabla {\cal T}^{(\beta)}(y+\teta)(1+k|y|^2)^{-n}\,dy=0\quad\hbox{if and only if}\quad \teta=0.\eqno\seiventisetteprimo$$
Then \seiventisetteprimo\ and \seiventisei\ yield $\xi=0$ and so from \seiventisette\ we get
$$\int_{\erre^n}{\cal T}^{(\beta)}(y)(1+k|y|^2)^{-n}\,dy=0.$$
On the other hand
$$\displaylines{\int_{\erre^n}{\cal T}^{(\beta)}(y)(1+k|y|^2)^{-n}\,dy=\sum_{j=1}^n\int_{\erre^n}a_j|y_j|^{\beta}(1+k|y|^2)^{-n}\,dy\cr
=\left(\sum_{j=1}^n a_j\right)\int_{\erre^n}|y_1|^{\beta}(1+k|y|^2)^{-n}\,dy\not=0\cr}$$
which gives a contradiction.\fine
\medskip\noindent
{\bf Remark 2.11.}\quad We can prove the same result of the previous theorem with the assumptions
$$\left\{\eqalign{&\widehat Q_i\in C^{[\beta]-1,1}(B_{\eps_0}(q_0)),\quad\widehat  Q_i(q_0)\geq {1\over A_2},\cr
                 &\widehat Q_i(y)=\widehat Q_i(0)+{\cal T}_i^{(\beta)}(y)+{\cal R}_i(y)\quad\hbox{for}\quad |y|\leq\eps_0,\cr}\right.$$
where $y$ is some geodesic normal coordinates centered at $q_0$,
$$\left\{\eqalign{&{\cal T}_i^{(\beta)}(\lambda y)=\lambda^\beta {\cal T}_i^{(\beta)}(y)\quad\forall\, \lambda>0,\quad y\in \erre^n,\cr
&\sum_{s=0}^{\scriptscriptstyle [\beta]}|\nabla^s{\cal R}_i(y)||y|^{-\beta+s}\mathop{\longrightarrow}\limits_{y\to0}0\quad\hbox{uniformly for}\ i,\cr}\right.$$
${\cal T}_i^{(\beta)}\to {\cal T}^{(\beta)}$ in $C^1(S^{n-1})$ and there exists a positive constant $A_3$ such that
$$\eqalign{&A_3|y|^{\beta-1}\leq|\nabla {\cal T}^{(\beta)}(y)|,\quad|y|\leq\eps_0,\cr
&\left(\matrix{\int_\erre^n\nabla {\cal T}^{(\beta)}(y\xi)(1+|y|^2)^{-n}\, dy\cr
       \int_\erre^n {\cal T}^{(\beta)}(y\xi)(1+|y|^2)^{-n}\, dy\cr}\right)\not=0\quad\forall\,\xi\in\erre^n.\cr}$$
\medskip\noindent
Theorems 2.8 and 2.9 easily imply the following result.
\medskip\noindent
{\bf Theorem 2.12.}\quad \sl For $n\geq 5$, suppose that $\widehat Q\in C^1(S^n)$ is some positive function such that for any critical point $q_0$ there exists some real number $\beta=\beta(q_0)$ with $\beta\geq 2$ and $\beta\in(n-4,n)$ such that in some geodesic normal coordinate system centered at $q_0$
$$\widehat Q(y)=\widehat Q(q_0)+\sum_{j=1}^n a_j|y_j|^\beta+R(y)\eqno\seiventotto$$
where $a_j=a_j(q_0)\not=0$, $\sum a_j\not=0$, $R(y)\in C^{[\beta]-1,1}$ near $0$ and
$$\sum_{s=0}^{[\beta]}|\nabla^s\widehat Q(y)||y|^{-\beta+s}=o(1)\quad\hbox{as}\ y\to 0.$$
Then for any $\eps>0$ there exists a positive constant $C(Q,n,\eps)$ such that for all $\eps\leq\mu\leq1$ any positive solution $v$ of $(P)$ with $\widehat Q$ replaced by 
$$\widehat Q_t=t\widehat Q+(1-t){(n^2-4)n\over 8}$$
satisfies 
$$C(Q,n,\eps)^{-1}\leq v(q)\leq c(Q,n,\eps),\quad\forall\, q\in S^n.$$
\medskip\noindent\rm
The following result is an existence theorem (which is the analogous of [\liuno,  Theorem 0.1]) which we will prove by making a homotopy which reduces the problem to the perturbative case.
\medskip\noindent
{\bf 2.3. A homotopy argument.}\quad We are able now to build a homotopy which allows us to pass from the perturbative case to the non-perturbative one.\medskip\noindent
{\bf Theorem 2.13.}\quad\sl For $n\geq 5$, suppose that $\widehat Q$ satisfies the assumptions of Theorem 2.12 and that
$$\sum_{\nabla_{g_0}\widehat Q(q_0)=0 \atop \sum_{j=1}^n a_j(q_0)<0}(-1)^{i(q_0)}\not=(-1)^n$$
where
$$i(q_0)=\#\{a_j(q_0):\ \nabla_{g_0}\widehat Q(q_0)=0,\ a_j(q_0)<0,\ 1\leq j\leq n\}.$$
Then (P) has at least one $C^4$ positive solution.
\medskip\noindent
{\bf Proof.}\rm \quad Let us consider $\widehat Q_t=t\widehat Q+(1-t)(n^2-4)n/8$ and the associated problems
$$\left\{\eqalign{&\pan v={n-4\over 2} \widehat Q_t v^{n+4\over n-4}\cr
                  &u>0\cr}\qquad \hbox{on}\ S^n.\right.\eqno(P_t)$$
Let $v_t$ be solutions of $(P_t)$, $0\leq t\leq1$. It follows from Theorem 2.8 that after passing to a subsequence either $\{v_t\}_{0\leq t\leq1}$ is bounded in $L^\infty(S^n)$ (and consequently in $C^{4,\alpha}(S^n)$, $0<\alpha<1$) or $\{v_t\}_{0\leq t\leq1}$ has precisely one isolated simple blow up point. Lemma 1.2 and Theorem 1.5 imply that for any $\delta>0$ we have solutions corresponding to $\widehat Q_\eps$ for $\eps$ small in the $\delta-$neighbourhood ($H^2(S^n)-$topology) of $Z_0$, where $Z_0$ is the set of the solutions of $(P_0)$
$$\pan v={(n-4)(n^2-4)n\over 16}v^{n+4\over n-4},\eqno(P_0)$$
i.e. up to a constant $Z_0$ is the critical manifold $Z$ of Section 1.\par
Let us fix $\eps>0$ and consider solutions of $(P_t)$ for $\eps\leq t\leq1$. Lemma 2.4 implies that solutions can blow up only at a precisely one of the critical points of $Q$. Theorem 2.12 implies that for $\eps\leq t\leq1$ any solution of $(P_t)$ $v$ satisfies
$$C^{-1}\leq v(x)\leq C\quad\forall\, x\in S^n$$
for some positive constant $C=C(Q,\eps,n)$, and then all the solutions of $(P_t)$ for $\eps\leq t\leq1$ are in a compact set of the space
$$C^{4,\alpha}(S^n)_+=\{w\in C^{4,\alpha}(S^n):\ w>0\ \hbox{on}\ S^n\}.$$
Then we can find some bounded open subset of $C^{4,\alpha}(S^n)_+$ denoted by ${\cal O}_\eps$ which contains all positive solutions of $(P_t)$ for $\eps\leq t\leq 1$. Thanks to the homotopy invariance of the Leray-Schauder degree we have
$$\eqalign{&\deg_{C^{4,\alpha}(S^n)_+}\left(v-(\pan)^{-1}\left({n-4\over 2}\widehat Qv^{n+4\over n-4}\right),{\cal O}_\eps,0\right)\cr
&\quad=\deg_{C^{4,\alpha}(S^n)_+}\left(v-(\pan)^{-1}\left({n-4\over 2}\widehat Q_\eps v^{n+4\over n-4}\right),{\cal O}_\eps,0\right).\cr}$$
Arguing as in [\cgydue, \dma, \liuno], one can prove that the right-hand side is equal to
$$\sum_{\nabla_{g_0}\widehat Q(q_0)=0 \atop \sum_{j=1}^n a_j(q_0)<0}(-1)^{i(q_0)}-(-1)^n$$
thus getting the conclusion.\fine
\medskip\noindent
{\bf Remark 2.14.}\quad \rm Let us point out the main differences between the above results and the results of [\dmadue]. In [\dmadue] $\widehat Q$ is supposed to be a Morse function satisfying a non-degeneracy condition (namely $\Delta_{g_0}\widehat Q(x)\not=0$ whenever $\nabla \widehat Q(x)=0$); no flatness condition is required. Under these assumptions when $n=5$ no blow up occurs. Instead the case $n=6$ can present multiple blow up points. Making some more assumptions on $\widehat Q$ and after a more precise description of the blow up scheme, an existence results when $n=6$ is stated in [\dmadue, Theorem 1.9]. In our case, the flatness assumption allows us to prove that there are no multiple blow up, thus getting compactness and allowing the homotopy argument for any $n\geq 5$.

\vskip2truecm
\nuovoparagrafo 3

\centerline{\bf 3.\quad The symmetric case}
\bigskip\noindent
\rm
We now consider the case in which the prescribed curvature is invariant under the action of a group $\Sigma\subset O(n+1)$. Let us start from the perturbative case. Suppose that $\widehat K$ is $\Sigma-$invariant. In this case we can argue as in [\simmetria] to find critical points of $f_\eps$ and as in Section 1 (Van der Vorst techniques) to prove the positivity of such critical points so that we get exactly the analogous results of Theorem 5.1, 5.3, 5.4, 6.1, and 6.4 of [\simmetria]. In particular we get
\medskip\noindent
{\bf Theorem 3.1.}\quad \sl Let $\fix(\Sigma)=\{p\in S^n:\ p_1=\dots=p_k=0\}$ for some $1\leq k\leq n$ and assume that $P_N$ is a nondegenerate minimum for $\widehat K$ on $S^n$. Moreover let $K=\widehat K\circ\pi$ satisfy
$$\eqalignno{\hbox{(i)}&\quad Y_k:=Crit\, (K)\cap V_{n-k}\ \hbox{is finite and all}\ \xi\in Y_k\ \hbox{are non-degenerate,}&\cr
\hbox{\ \ \ }&\quad \hbox{where}\ V_{n-k}=\{\xi\in\erre^n:\ \xi_1=\dots=\xi_k=0\};&\cr
\hbox{(ii)}&\quad \Delta K(\xi)\not=0\quad\forall\,\xi\in Y_k\quad\hbox{and}&\cr
\hbox{\ \ \ }&\quad\qquad\quad\quad \sum_{\xi\in Y_k,\ \Delta K(\xi)<0}(-1)^{\tilde m}-(-1)^{n-k}\not=0&\grado\cr
\hbox{\ \ \ }&\quad\hbox{where} \ \tilde m\ \hbox{is the Morse index of}\ \xi\ \hbox{as a critical point of}\ K\ \hbox{constrained} &\cr
\hbox{\ \ \ }&\quad \hbox{on}\ V_{n-k}.&\cr}$$
Then \tre\ has a symmetric solution.\fine\rm\medskip\noindent
{\bf Remark 3.2.}\quad In an analogous way, using Lemma 1.2 and the arguments of [\simmetria], we can prove that, if $\widehat Q$ is of the form \seiventotto\ near the critical fixed points, \tre\ has a symmetric solution provided
$$\sum_{q_0\in {\rm Crit}\,(Q)\cap \fix (\Sigma)\atop \sum a_j(q_0)<0}(-1)^{i_\Sigma(q_0)}-(-1)^{n-k}\not=0$$
where $i_\Sigma(q_0)$ denotes the Morse index of $q_0$ as a critical point of $Q$ restricted to $\fix(\Sigma)$.
\medskip\noindent
We can treat the non perturbative case as Ambrosetti-Li-Malchiodi did in [\alm] for the scalar curvature problem, thus getting the following
\medskip\noindent
{\bf Theorem 3.3.}\quad \sl For $n\geq 5$, let $\widehat Q\in C^{n-4,\alpha}(S^n)$ positive and invariant under the action of a group $\Sigma \subseteq O(n+1)$, such that
$$\leqalignno{&\hbox{for all critical points of}\ \widehat Q \ \hbox{in}\ \fix(\Sigma):=\{q\in S^n:\ \sigma q=q\quad \forall\,\sigma \in \Sigma\}&(H1)\cr
& \widehat Q\ \hbox {is of the form \seiventotto;}&\cr
&\widehat Q\ \hbox {satisfies}\  (*)_\beta \ \hbox{near the critical points;}&(H2)\cr}$$
suppose that 
$$\eqalignno{&{\rm Crit}\,(\widehat Q)\cap \fix(\Sigma)\quad\hbox{is finite,}&\seiventinove\cr
             &\fix(\Sigma)=\{x=(x_1,\dots,x_{n+1})\in S^n:\ x_1=\dots=x_k=0\}.&\seitrentuno\cr}$$
Then all $\Sigma-$invariant solutions of $(P)$ stay in a compact set of 
$$C^{4,\alpha}(S^n)^+_{\Sigma}=\{w\in C^{4,\alpha}(S^n)^+: w\ \hbox{is}\ \Sigma-\hbox{invariant}\}$$
and 
$$\eqalignno{&\deg_{C^{4,\alpha}(S^n)^+_{\Sigma}}\left(v-{\pan}^{-1}\left({n-4\over 2}\widehat Qv^{n+4\over n-4},{\cal O},0\right)\right)&\cr
&\qquad\qquad =\sum_{q_0\in {\rm Crit}\,(Q)\cap \fix (\Sigma)\atop \sum a_j(q_0)<0}(-1)^{i_\Sigma(q_0)}-(-1)^{n-k}&\seitrentadue\cr}$$
where ${\cal O}$ is some bounded open set in $C^{4,\alpha}(S^n)^+_{\Sigma}$ containing all $\Sigma-$invariant solutions of $(P)$. In particular, if the number on the right hand side of \seitrentadue\ is nonzero, $(P)$ has at least one $\Sigma-$invariant solution.
\medskip\noindent\rm
{\bf Proof.}\quad We can repeat exactly the first part of the proof of the previous theorem. Since there is at most one blow up point, we know that $\Sigma-$invariant solutions of $(P_t)$ for $0\leq t\leq1$ can not blow up at any point in $S^n\,\backslash\, \fix(\Sigma)$. \par
Fix now $\eps>0$ and consider $(P_t)$ for $\eps\leq t\leq1$. Away from $\fix(\Sigma)$, $\Sigma-$invariant solutions stay uniformly bounded. Lemma 2.4 implies that solutions can blow up only at precisely one of the critical points of $Q$ in $\fix(\Sigma)$. Theorem 2.12 implies that there exists a bounded open subset of $C^{4,\alpha}(S^n)^+_{\Sigma}$ denoted as ${\cal O}_\eps$ which contains  all $\Sigma-$invariant solutions of $(P_t)$ for $\eps\leq t\leq1$. As above, the homotopy invariance of the Leray-Schauder degree allows us to conclude the proof.\fine

\vskip2truecm
\nuovoparagrafo 4

\centerline{\bf 4.\quad Appendix}
\bigskip\noindent
Let us give the proof of the technical lemmas stated in Section 2.
\medskip\noindent 
\rm{\bf Proof of Lemma 2.4.}\quad Let us take a cut-off function $\eta\in C^{\infty}_c(B_1)$ satisfying
$$\eta(x)=\left\{\eqalign{1,\quad&\hbox{if}\ |x|\leq\quarto\cr
                           0,\quad&\hbox{if}\ |x|\geq\mez\cr}\right.$$
and multiply $\seidue_i$ by $\eta{\partial \ui\over\partial x_j}$
$$\int_{B_1}\Delta^2\ui\eta{\partial\ui\over\partial x_j}={n-4\over 2}\int_{B_1}Q_i\ui^p{\partial\ui\over\partial x_j}\eta.$$
Integrating by parts on $B_1$ we get
$$\int_{B_1}Q_i\ui^p{\partial\ui\over\partial x_j}\eta=-{1\over p+1}\int_{B_1}{\partial Q_i\over\partial x_j}\eta\ui^{p+1}-{1\over p+1}\int_{B_1}Q_i{\partial\eta\over\partial x_j}\ui^{p+1}$$
and
$$\eqalign{&\int_{B_1}\Delta^2\ui\eta{\partial\ui\over\partial x_j}\cr
&\quad=\int_{B_1}\Delta\ui\Delta\eta{\partial\ui\over\partial x_j}+\int_{B_1}\Delta\ui\eta\Delta\left({\partial \ui\over\partial x_j}\right)+2\int_{B_1}\Delta\ui\left\langle\nabla\eta,\nabla\left({\partial \ui\over\partial x_j}\right)\right\rangle\cr
&\quad=\int_{B_1}\Delta\ui\Delta\eta{\partial\ui\over\partial x_j}-\mez\int_{B_1}(\Delta\ui)^2{\partial\eta\over\partial x_j}+2\int_{B_1}\Delta\ui\left\langle\nabla\eta,\nabla\left({\partial \ui\over\partial x_j}\right)\right\rangle.\cr}$$
Recalling that $\nabla \eta$ and $\Delta\eta$ vanish in $B_1\,\backslash\,\big(B_{\mez}\,\backslash\, B_{\quarto}\big)$ the two previous equalities yield
$$\eqalignno{&{1\over \pii+1}{n-4\over 2}\int_{B_{\mez}}{\partial Q_i\over\partial x_j}\ui^{\pii+1}\eta=-\int_{B_{\mez}\,\backslash\, B_{\quarto}}\Delta\ui{\partial\ui\over\partial x_j}\Delta\eta+\mez\int_{B_{\mez}\,\backslash\, B_{\quarto}}{\partial \eta\over\partial x_j}(\Delta\ui)^2&\cr
&\quad-2\int_{B_{\mez}\,\backslash\, B_{\quarto}}\Delta\ui\left\langle\nabla\eta,\nabla\left({\partial\ui\over\partial x_j}\right)\right\rangle-{n-4\over 2}{1\over \pii+1}\int_{B_{\mez}\,\backslash\, B_{\quarto}}{\partial \eta\over\partial x_j}Q_i\ui^{\pii+1}.&\seinove\cr}$$
Owing to the weak Harnack inequality proved in [\dmadue, Lemma 2.5] and to Schauder's elliptic estimates and using \propduenove, from \seinove\ we get
$$\left|\int_{B_{\mez}}{\partial Q_i\over\partial x_j}\ui^{\pii+1}\right|\leq c_1\ui^{-2}(\yi)+c_2\ui^{-\pii-1}(\yi)\leq c_3\ui^{-2}(\yi).\eqno\seiundici$$
From \seiundici, it follows
$$\left|{\partial Q_i\over\partial x_j}(\yi)\int_{B_{\mez}}\ui^{\pii+1}\right|-c_1\ui^{-2}(\yi)\leq c_2\left|\int_{B_{\mez}}\ui^{\pii+1}\left\{{\partial Q_i\over\partial x_j}(\yi)-{\partial Q_i\over\partial x_j}(y)\right\}\right|.$$
Just to simplify notations, let us suppose that $\beta\in\enne$; it is clear the all the following arguments can be played also in the general case. Expanding $\nabla Q_i$ and using (Q2) we deduce
$$\eqalign{&\left|{\partial Q_i\over\partial x_j}(\yi)\int_{B_{\mez}}\ui^{\pii+1}\right|-c_1\ui^{-2}(\yi)\cr
&\quad \leq c_2\int_{B_{\mez}}\Big\{\sum_{s=2}^{\scriptscriptstyle\beta-1}|\nabla^sQ_i(\yi)||y-\yi|^{s-1}+\max_{0\leq t\leq1}|\nabla^{\beta}Q_i(\yi+t(y-\yi))||y-\yi|^{\beta-1}\Big\}\ui^{\pii+1}\cr
&\quad\leq c_3\int_{B_{\mez}}\Big\{\sum_{s=2}^{\beta-1}|\nabla Q_i(\yi)|^{\beta-s\over \beta-1}|y-\yi|^{s-1}+c_4|y-\yi|^{\beta-1}\Big\}\ui^{\pii+1}.\cr
}$$
Hence
$$|\nabla Q_i(\yi)|\leq c_1\ui(y_i)^{-2}+c_2\int_{B_{\mez}}\Big\{\sum_{s=2}^{\scriptscriptstyle\beta-1}|\nabla Q_i(\yi)|^{\beta-s\over \beta-1}|y-\yi|^{s-1}+c_3|y-\yi|^{\beta-1}\Big\}\ui^{\pii+1}.$$
By the Young inequality with conjugate exponents $q={\beta-1\over\beta- s}$ and $q'={\beta-1\over s-1}$ we have that for any $\delta>0$
$$|\nabla Q_i(\yi)|^{\beta-s\over \beta-1}|y-\yi|^{s-1}\leq {\beta-s\over \beta-1}\delta^{(\beta-1)/( \beta-s)}|\nabla Q_i(\yi)|+{s-1\over\beta-1}\left({1\over\delta}\right)^{\beta-1\over s-1}|y-\yi|^{\beta-1}$$
and hence
$$\eqalignno{|\nabla Q_i(\yi)|\leq &\,c_1\ui^{-2}(\yi)+c_2\delta^{\beta-1\over \beta-s}|\nabla Q_i(\yi)|\int_{B_{\mez}}\ui^{p+1}&\cr
&+c_4\delta^{-{\beta-1\over s-1}}\int_{B_{\mez}}|y-\yi|^{\beta-1}\ui^{p+1}+c_5\int_{B_{\mez}}|y-\yi|^{\beta-1}\ui^{\pii+1}.&\seidodici\cr}$$
Choosing $\delta$ small enough and recalling that $\stima_s$ with $s=0$ yields $\int_{B_{\mez}}\ui^{p+1}\leq {\rm const}$ from \seidodici\ and $\stima_s$ with $s=\beta-1$ again we deduce 
$$\eqalignno{|\nabla Q_i(\yi)|&\leq c_1\ui^{-2}(\yi)+\mez|\nabla Q_i(\yi)|+c_2\int_{B_{\mez}}|y-\yi|^{\beta-1}\ui^{\pii+1}&\cr
& \leq c_1\ui^{-2}(\yi)+\mez|\nabla Q_i(\yi)|+c_3\ui^{-{2\over n-4}(\beta-1)}(\yi).&\seitredici\cr}$$
Estimate \seitredici\ yields the conclusion of the lemma for $\beta\in\enne$. In the general case the estimate can be proved with analogous calculations. \fine
\medskip\noindent
{\bf Proof of Lemma 2.5.}\quad\rm We have that
$$\eqalignno{\left|\int_{B_\sigma(\yi)}(y-\yi)\cdot\nabla Q_i\ui^{\pii+1}\right|&\leq \left|\nabla Q_i(\yi)\cdot\int_{B_\sigma(\yi)}(y-\yi)\ui^{\pii+1}\right|&\cr
&\quad +\int_{B_\sigma(\yi)}|y-\yi||\nabla (Q_i)(y)-\nabla Q_i(\yi)|\ui^{\pii+1}.&\seiquattordici\cr}$$
The change of variable $y=\ui^{-{p-1\over 4}}(\yi)x+\yi$ yields 
$$\int_{B_\sigma(\yi)}(y-\yi)\ui^{\pii+1}=\ui^{-{(p-1)(n+1)\over 4}}(\yi)\int_{B_{a_i}}x\ui^{p+1}\Big(\ui^{-{p-1\over 4}}(\yi)x+\yi\Big)\,dx\eqno\cambiouno$$
where $a_i={\sigma\ui^{(\pii-1)/4}}(\yi)$. Using [\dmadue, Proposition 2.7] we have that, up to a subsequence,
$$J_i(x):={\ui^{p+1}\Big(\ui^{-{p-1\over 4}}(y_i)x+\yi\Big)\over\ui^{-p-1}(\yi)}-\left({1\over 1+k_i|x|^2}\right)^n,\qquad k_i^2=[2n(n-2)(n+2)]^{-1}Q_i(\yi),$$
satisfies
$$\|J_i\|_{C^2(B_{2R_i})}\leq\eps_i\eqno\propduesette$$
for any $R_i\to\infty$ and $\eps_i\to 0^+$. Then \cambiouno\ gives
$$\eqalignno{\left|\int_{B_\sigma(\yi)}(y-\yi)\ui^{\pii+1}\right|&=\left|\ui^{-{2\over n-4}}(y_i)\int_{B_{a_i}}x[(1+k_i|x|^2)^{-n}+J_i(x)]\, dx\right|\cr
               &=\left|\ui^{-{2\over n-4}}(\yi)\int_{B_{a_i}}xJ_i(x)\,dx\right|\leq\ui^{-{2\over n-4}}(\yi)\eps_ia_i|B_{a_i}|.\cr}$$
Since $\{\eps_i\}_i$ is an arbitrary sequence going to $0$ we get that
$$\nabla Q_i(\yi)\cdot\int_{B_\sigma(\yi)}(y-\yi)\ui^{\pii+1}=o\left(|\nabla Q_i(\yi)|\ui^{-{2\over n-4}}(\yi)\right).\eqno\cambiodue$$
Suppose for simplicity $\beta\in\enne$. Using (Q2) and expanding $\nabla Q_i$ in a neighbourhood of $y_i$, we get
$$\eqalign{&\left|{\partial Q_i\over\partial y_j}(y)-{\partial Q_i\over\partial y_j}(\yi)\right|\cr
&\quad\leq c_1\left\{\sum_{s=2}^{\scriptscriptstyle\beta-1}|\nabla^s Q_i(\yi)||y-\yi|^{s-1}+\max_{0\leq t\leq 1}|\nabla^{\beta}Q_i(\yi+t(y-\yi))|\cdot|y-\yi|^{\beta-1}\right\}\cr
&\quad\leq c_2\left\{\sum_{s=2}^{\scriptscriptstyle\beta-1}|\nabla Q_i(\yi)|^{\beta-s\over\beta-1}|y-\yi|^{s-1}+|y-\yi|^{\beta-1}\right\}\cr
&\quad \leq c_3\sum_{s=2}^{\scriptscriptstyle\beta}|\nabla Q_i(\yi)|^{\beta-s\over\beta-1}|y-\yi|^{s-1}\cr}$$
and so
$$\int_{B_\sigma(\yi)}|y-\yi|\cdot\left|{\partial Q_i\over\partial y_j}(y)-{\partial Q_i\over\partial y_j}(\yi)\right|\ui^{\pii+1}\leq c\int_{B_\sigma(\yi)}\sum_{s=2}^{\scriptscriptstyle\beta}|\nabla Q_i(\yi)|^{\beta-s\over\beta-1}|y-\yi|^{s}\ui^{p+1}.\eqno\seiotto$$
Note that by means of the Young inequality with the conjugate exponents $q={\beta-1\over\beta-s}$ and $q'={\beta-1\over s-1}$ we get
$$\eqalignno{|\nabla Q_i(\yi)|^{\beta-s\over\beta-1}|y-\yi|^{s}&=|y-\yi||\nabla Q_i(\yi)|^{\beta-s\over\beta-1}|y-\yi|^{s-1}&\cr
&\leq {\beta-s\over\beta-1}|y-\yi||\nabla Q_i(\yi)|+{s-1\over\beta-1}|y-\yi|^\beta.&\seiottoprimo\cr}$$
Then \seiotto\ and \seiottoprimo\ yield
$$\eqalignno{&\int_{B_\sigma(\yi)}|y-\yi|\cdot\left|{\partial Q_i\over\partial y_j}(y)-{\partial Q_i\over\partial y_j}(\yi)\right|\ui^{\pii+1}&\cr
&\qquad\leq c_1|\nabla Q_i(\yi)|
\int_{B_\sigma(\yi)}|y-\yi|\ui^{\pii+1}+c_2\int_{B_\sigma(\yi)}|y-\yi|^\beta\ui^{\pii+1}.&\seiottosecondo\cr}$$
If $\beta=2$, from \seiquattordici, \cambiodue, \seiotto, and $\stima_s$ with $s=2$ we get
$$\eqalign{\left|\int_{B_\sigma(\yi)}(y-\yi)\cdot\nabla Q_i\ui^{\pii+1}\right|&\leq o\left(|\nabla Q_i(\yi)|\ui^{-{2\over n-4}}(\yi)\right)+c_1\int_{B_\sigma(\yi)}|y-\yi|^2\ui^{\pii+1}\cr
&\leq o\left(|\nabla Q_i(\yi)|\ui^{-{2\over n-4}}(\yi)\right)+c_2\ui^{-{4\over n-4}}(\yi).\cr}$$
Consider now the case $\beta>2$: applying $\stima_s$ with $s=1$ and $s=\beta$ and using \seiquattordici, \cambiodue, and \seiottosecondo\ we get
$$\eqalign{&\left|\int_{B_\sigma(\yi)}(y-\yi)\nabla Q_i\ui^{\pii+1}\right|\cr
&\quad\leq c_1|\nabla Q_i(\yi)|\ui^{-{2\over n-4}}(\yi)+c_2\int_{B_\sigma(\yi)}|y-\yi||\nabla Q_i(y)-\nabla Q_i(\yi)|\ui^{\pii+1}\cr
&\quad \leq c_1|\nabla Q_i(\yi)|\ui^{-{2\over n-4}}(\yi)+c_3|\nabla Q_i(\yi)|\int_{B_\sigma(\yi)}|y-\yi|\ui^{\pii+1}+c_4\int_{B_\sigma(\yi)}|y-\yi|^{\beta}\ui^{\pii+1}\cr
&\quad \leq c_5|\nabla Q_i(\yi)|\ui^{-{2\over n-4}}(\yi)+c_6\ui^{-{2\beta\over n-4}}(\yi).\cr}$$
As before, the above arguments can be played also in the general case, i.e. for $\beta\not\in\enne$.\fine
\medskip\noindent\rm

\vskip2truecm
\centerline {\bf References}\bigskip
\item {[\pertuno]} \maiuscolo A. Ambrosetti and M. Badiale: \sl Homoclinics: Poincar\'e-Melnikov type results via a variational approach, \rm Ann. Inst. H. Poincar\'e Anal. Non Lin\'eaire, \bf 15 \rm(1998), 233-252.
\smallskip
\item {[\pertdue]} \maiuscolo A. Ambrosetti and M. Badiale: \sl Variational perturbative methods and bifurcation of bound states from the essential spectrum, \rm Proc. Royal Soc. Edinburgh, \bf  128A\rm(1998), 1131-1161.
\smallskip
\item {[\spagnoli]} \maiuscolo A. Ambrosetti, J. Garcia Azorero, and I. Peral: \sl Perturbation of $\Delta u+u^{N+2\over N-2}=0$ the scalar curvature problem in $\erre^N$ and related topics, \rm J. Funct. Analysis, \bf  165\rm(1999), 117-149.
\smallskip
\item {[\alm]} \maiuscolo A. Ambrosetti, Y. Y. Li, and A. Malchiodi: \sl A note on the scalar curvature problem in the presence of symmetries, \rm Ricerche di Mat., to appear.
\smallskip
\item {[\simmetria]} \maiuscolo A. Ambrosetti and A. Malchiodi: \sl On the symmetric scalar curvature problem on $S^n$, \rm J. Diff. Equations, \bf 170 \rm(2001), 228-245.
\smallskip
\item {[\aubin]} \maiuscolo T. Aubin: \rm ``Some nonlinear problems in differential geometry'', Springer-Verlag, 1998.
\smallskip
\item {[\bgm]} \maiuscolo M. Berger, P. Gauduchon, and E. Mazet: \rm ``Le spectre d'une vari\'et\'e riemannienne'', Lecture Note in Mathematics, \bf 194, \rm Springer-Verlag, New York/Berlin, 1971.
\smallskip
\item {[\branson]} \maiuscolo T.P. Branson: \sl Differential operators canonically associated to a conformal structure, \rm Mathematica Scandinavica, \bf 57 \rm(1985), .
\smallskip
\item {[\bcy]} \maiuscolo T.P. Branson, S.A. Chang, and P.C. Yang: \sl Estimates and extremal problems for the log-determinant on 4-manifolds, \rm Comm. Math. Phys. \bf 149 \rm(1992), 241-262.
\smallskip
\item {[\chang]} \maiuscolo S.A. Chang: \sl On Paneitz operator-a fourth order differential operator in conformal geometry, \rm Survey article, to appear in the Proceedings for the 70th birthday of A.P. Calderon.
\smallskip
\item {[\cgydue]} \maiuscolo S.A. Chang, M.J. Gursky, and P.C. Yang: \sl The scalar curvature equation on $2-$ and $3-$sphere, \rm Calc. Var. Partial Differential Equations, \bf 1 \rm(1993), 205-229.
\smallskip
\item {[\cgy]} \maiuscolo S.A. Chang, M.J. Gursky, and P.C. Yang: \sl Regularity of a fourth order non linear PDE with critical exponent, \rm Amer. J. Math., \bf 121-1 \rm(1999), 215-257.
\smallskip
\item {[\cqy]} \maiuscolo S. A. Chang, J. Qing, and P. Yang: \sl On the Chern-Gauss-Bonnet integral for conformal metrics on $\erre^4$, \rm Duke Math., to appear.
\smallskip
\item {[\changyangdue]} \maiuscolo S.A. Chang and P.C. Yang: \sl A perturbation result in prescribing scalar curvature on $S^n$, \rm Duke Mathematical Journal, \bf 64\rm (1991), 27-69.
\smallskip
\item {[\changyang]} \maiuscolo S.A. Chang and P.C. Yang: \sl On a fourth order curvature invariant, \rm Contemporary Mathematics, \bf 237, \rm Spectral problems in Geometry and Arithmetic, Ed. T. Branson, AMS, 1999, 9-28.
\smallskip
\item {[\dhl]} \maiuscolo Z. Djadli, E. Hebey, and M. Ledoux: \sl Paneitz type operators and applications, \rm Duke Mathematical Journal, \bf 104\rm(2000), no. 1, 129-169.
\smallskip
\item {[\dma]} \maiuscolo Z. Djadli, A. Malchiodi, and M. Ould Ahmedou: \sl Prescribed fourth order conformal invariant on the standard sphere, part I: a perturbative result, \rm Preprint.
\smallskip
\item {[\dmadue]} \maiuscolo Z. Djadli, A. Malchiodi, and M. Ould Ahmedou: \sl Prescribed fourth order conformal invariant on the standard sphere, part II: blow up analysis and applications, \rm Pre\-print.
\smallskip
\item {[\gursky]} \maiuscolo M.J. Gursky: \sl The Weyl functional, de Rham cohomology, and Kahler-Einstein metrics, \rm Ann. of Math., \bf 148 \rm(1998), 315-337.
\smallskip  
\item {[\hebey]}
\maiuscolo E. Hebey: \rm ``Introduction \`a l'analyse non lin\'eaire sur les vari\'et\'es'', Diderot Editeur, Paris, 1997. 
\smallskip 
 \item {[\liuno]} \maiuscolo Y. Y. Li,
\sl Prescribing scalar curvature on $S^n$ and related topics, Part I, \rm J.
Differential Equations, \bf  120\rm(1995), 319-410.
 \smallskip 
 \item {[\lidue]} \maiuscolo Y. Y. Li, \sl Prescribing scalar curvature on $S^n$ and
related topics, Part II: Existence and compactness, \rm Comm. Pure Appl.
Math., \bf  49\rm(1996), 437-477.
\smallskip
\item {[\lincs]} \maiuscolo C. S. Lin: \sl A classification of solutions of a conformally invariant fourth order equation in $\erre^n$, \rm Comment. Math. Helv., \bf 73 \rm(1998), 206-231.
 \smallskip
\item {[\paneitz]} \maiuscolo S. Paneitz: \sl A quartic conformally covariant differential operator for arbitrary pseu\-do\--Riemannian manifolds, \rm Preprint, 1983.
\smallskip
\item {[\van]} \maiuscolo R. C. A. M. Van Der Vorst: \sl Best constants for the embedding of the space $H^2\cap H^1_0(\Omega)$ into $L^{2N\over (N-4)}(\Omega)$, \rm Differential and Integral Equations, \bf 6 \rm(1993), 259-276.
\smallskip
\item {[\weixu]} \maiuscolo J. Wei and X. Xu: \sl On conformal deformations of metrics on $S^n$, \rm Journal of Functional Analysis, \bf 157 \rm(1998), no. 1, 292-325.
\smallskip
\end